\documentclass[journal]{IEEEtran}

\usepackage[utf8]{inputenc}
\usepackage{color}
\usepackage{array}
\usepackage{verbatim}
\usepackage{float}
\usepackage{amsmath}
\usepackage{amsthm}
\usepackage{amssymb}
\usepackage{graphicx}
\usepackage{longtable}
\usepackage{multirow}
\usepackage[unicode=true,
 bookmarks=false,
 breaklinks=false,pdfborder={0 0 1},colorlinks=false]
 {hyperref}
\hypersetup{
 colorlinks,bookmarksopen,bookmarksnumbered,citecolor=blue,urlcolor=blue}
\usepackage{cite}
\makeatletter
 \let\oldforeign@language\foreign@language
 \DeclareRobustCommand{\foreign@language}[1]{%
   \lowercase{\oldforeign@language{#1}}}

 \let\oldforeign@language\foreign@language
 \DeclareRobustCommand{\foreign@language}[1]{%
   \lowercase{\oldforeign@language{#1}}}

\ifCLASSINFOpdf
\else
\fi

\newtheorem{lem}{Lemma}

\newtheorem{prop}{Proposition}
\newtheorem{thm}{Theorem}
\newtheorem{rem}{Remark}

\newtheorem{assum}{Assumption}

\begin{document}
	\bstctlcite{IEEEexample:BSTcontrol}
%
%
\onecolumn
\noindent\rule{18.1cm}{2pt}\\
\underline{To cite this article:}
{\bf{\textcolor{red}{H. A. Hashim, L. J. Brown, and K. McIsaac, "Guaranteed Performance of Nonlinear Attitude Filters on the Special Orthogonal Group SO(3)," IEEE Access, vol. 7, no. 1, pp. 3731–3745, December 2019.}}}\\
\noindent\rule{18.1cm}{2pt}\\

\noindent{\bf The published version (DOI) can be found at:  \href{https://doi.org/10.1109/ACCESS.2018.2889612}{10.1109/ACCESS.2018.2889612} }\\

\vspace{40pt}\noindent Please note that where the full-text provided is the Author Accepted Manuscript or Post-Print version this may differ from the final Published version. { \bf To cite this publication, please use the final published version.}\\

\textbf{
	\begin{center}
		Personal use of this material is permitted. Permission from the author(s) and/or copyright holder(s), must be obtained for all other uses, in any current or future media, including reprinting or republishing this material for advertising or promotional purposes.\vspace{420pt}\\
	\end{center}
 }
\small { \bf
	\vspace{30pt}\noindent Please contact us and provide details if you believe this document breaches copyrights. We will remove access to the work immediately and investigate your claim.
} 
\normalsize

\twocolumn
\title{Guaranteed Performance of Nonlinear Attitude Filters on the Special Orthogonal Group SO(3)}

\author{Hashim~A.~Hashim, Lyndon J. Brown, and~Kenneth McIsaac
\thanks{This work was supported in part by the Canadian Space Agency FAST program.}
\thanks{H. A. Hashim, L. J. Brown and K. McIsaac are with the Department of Electrical and Computer Engineering,
University of Western Ontario, London, ON, Canada, N6A-5B9, e-mail: hmoham33@uwo.ca, lbrown@uwo.ca and kmcisaac@uwo.ca.}
}


\markboth{--,~Vol.~-, No.~-, \today}{Hashim \MakeLowercase{\textit{et al.}}:Guaranteed Performance of Nonlinear Attitude Filters on the Special Orthogonal Group SO(3)}
\markboth{}{Hashim \MakeLowercase{\textit{et al.}}: Guaranteed Performance of Nonlinear Attitude Filters on the Special Orthogonal Group SO(3)}

\maketitle

\begin{abstract}
This paper proposes two novel nonlinear attitude filters evolved directly
on the Special Orthogonal Group $\mathbb{SO}\left(3\right)$ able
to ensure prescribed measures of transient and steady-state performance.
The tracking performance of the normalized Euclidean distance of attitude
error is trapped to initially start within a large set and converge
systematically and asymptotically to the origin from almost any initial
condition. The convergence rate is guaranteed to be less than the
prescribed value and the steady-state error does not exceed a predefined
small value.%
{} The first filter uses a set of vectorial measurements with the need
for attitude reconstruction. The second filter does not require attitude
reconstruction and instead uses only a rate gyroscope measurement
and two or more vectorial measurements. These filters provide good
attitude estimates with superior convergence properties and can be
applied to measurements obtained from low cost inertial measurement
units (IMUs). Simulation results illustrate the robustness and effectiveness
of the proposed attitude filters with guaranteed performance considering
high level of uncertainty in angular velocity along with body-frame
vector measurements.
\end{abstract}

\begin{IEEEkeywords}
Nonlinear complementary filter, Attitude estimator, observer, estimates, special orthogonal group,  error function,
prescribed performance function, systematic convergence, transformed error, steady-state error,
transient error, SO(3), PPF, IMUs.

\end{IEEEkeywords}

\IEEEpeerreviewmaketitle{}

\section{Introduction}

%
%
%
%
\IEEEPARstart{A}{ttitude} estimation of rigid-body systems plays an essential role
in many engineering applications such as robotics, aerial and underwater
vehicles and satellites. The orientation of the rigid-body can be
reconstructed algebraically given that two or more known inertial
vectors as well as their body-frame vectors are available at each
time instant for measurement, for example using TRIAD or QUEST algorithms
\cite{black1964passive,shuster1981three} and singular value decomposition
(SVD) \cite{markley1988attitude}. Nonetheless, body-frame measurements
are corrupted with unknown constant bias and random noise components
and the static estimation in \cite{black1964passive,shuster1981three,markley1988attitude}
provides unsatisfactory results, in particular, if the moving body
is equipped with low-cost inertial measurement units (IMUs) \cite{crassidis2007survey,mahony2008nonlinear}.

During the last few decades, a remarkable effort has been done to
achieve higher filtering performance with noise reduction through
Gaussian filters. One of the earliest detailed derivations of a Gaussian
filter is the extended Kalman filter (EKF) in \cite{lefferts1982kalman}.
A novel Kalman filter was proposed later in \cite{choukroun2006novel}
and showed better results in comparison with the EKF in \cite{lefferts1982kalman}.
Also, other Gaussian filters have been proposed, such as multiplicative
EKF (MEKF) \cite{markley2003attitude,stovner2018attitude}, invariant
EKF \cite{bonnable2009invariant}, and geometric approximate minimum-energy
filter \cite{zamani2013minimum}. A good survey of Gaussian attitude
filters can be found in \cite{crassidis2007survey,hashim2018SO3Stochastic}. However, nonlinear
deterministic attitude filters have better tracking performance, and
require less computational power when compared to Gaussian filters
\cite{crassidis2007survey,mahony2008nonlinear,hashim2018SO3Stochastic}. Accordingly, nonlinear
deterministic attitude filters received considerable attention \cite{crassidis2007survey,mahony2008nonlinear,hashim2018SO3Stochastic}.

The need for attitude filters robust against uncertainty in measurement
sensors, especially with the development of low-cost IMUs, played
a significant role in the development of nonlinear attitude filters,
for example \cite{mahony2008nonlinear,hamel2006attitude,zlotnik2017nonlinear,grip2012attitude,hashim2018SO3Stochastic,hashim2018Conf1,wu2018fast}.
These filters can be easily fitted knowing a rate gyroscope measurement
and two or more vectorial measurements taken, for instance, by low-cost
IMUs. In general, the nonlinear attitude filter is achieved via careful
selection of the error function. The selected error function in \cite{mahony2005complementary}
underwent slight modifications in \cite{mahony2008nonlinear,hamel2006attitude,grip2012attitude},
overall performance, however, was not significantly changed. The main
problem of the error function in \cite{mahony2005complementary,mahony2008nonlinear,hamel2006attitude,grip2012attitude}
consists in the slow convergence, especially with large initial attitude
error. A new form of the error function presented in \cite{zlotnik2017nonlinear,lee2012exponential}
offered faster error convergence to the equilibrium point. In addition, recently proposed robust nonlinear stochastic attitude filters offer fast convergence of the attitude error to small neighborhood of the equilibrium point \cite{hashim2018SO3Stochastic,hashim2018Conf1}. However, no systematic
convergence is observed in \cite{zlotnik2017nonlinear,lee2012exponential,hashim2018SO3Stochastic,hashim2018Conf1}
in other words, the transient performance does not follow a predefined
trajectory and the steady-state error can not be controlled. Thus,
the prediction of transient and steady-state error performance is
almost impossible.

Prescribed performance signifies trapping the error to initiate arbitrarily
within a given large set and reduce systematically and smoothly to
a given small residual set \cite{bechlioulis2008robust}. The convergence
of the error is constrained by a specified range during the transient
as well as the steady-state performance. The aim of prescribed performance
is to relax the constrained error and transform it to a new unconstrained
form. Accordingly, the new form allows one to keep the error below
the predefined value which could be useful in the estimation and control
process. Prescribed performance has been implemented successfully
in many control applications such as two degree of freedom planar
robot \cite{bechlioulis2008robust,mohamed2014improved}, nonlinear
control with input saturation \cite{zhu2018robust}, and uncertain
multi-agent system \cite{hashim2017neuro,hashim2017adaptive}. Attitude
error function is an essential step for the construction of any nonlinear
attitude filter, as it is directly related to the convergence behavior
of the error trajectory.

Accordingly, two robust nonlinear attitude filters on the Special Orthogonal
Group $\mathbb{SO}\left(3\right)$ with predefined transient as well
as steady-state characteristics are proposed in this paper. An alternate
attitude error function is selected such that the error is defined
in terms of normalized Euclidean distance. The error function is forced
to be contained and start within a predefined large set and reduce
systematically and smoothly to a known small set. Therefore, the aforementioned
error is constrained and as it approaches zero the transformed error,
which is a new form of unconstrained error, approaches the origin
and vice versa. These filters ensure boundedness of the closed loop
error signals with attitude error being regulated asymptotically to
the origin. The attitude estimators ensure faster convergence properties
and satisfy prescribed performance better than similar estimators
considered in the literature. The fast convergence is mainly attributed
to the behavior of the estimator gains, which are dynamic. The first
filter needs a rate gyroscope measurement and a set of two or more
vectorial measurements to obtain online algebraic reconstruction of
the attitude. The second filter uses the rate gyroscope measurement
combined with the aforementioned vectorial measurements directly avoiding
the need for attitude reconstruction.

The remainder of the paper is organized as follows: Section \ref{sec:SO3PPF_Math-Notations}
gives a brief review of the mathematical notation, $\mathbb{SO}\left(3\right)$
parameterization, and a number of selected relevant identities. Section
\ref{sec:SO3PPF_Problem-Formulation-in} formulates the attitude problem,
presents the estimator structure and error criteria, and formulates
the attitude error in terms of prescribed performance. The two proposed
filters and the associated stability analysis are demonstrated in
Section \ref{sec:SO3PPF-Filters}. Section \ref{sec:SO3PPF_Simulations}
illustrates through simulation the effectiveness and robustness of
the proposed filters. Finally, Section \ref{sec:SO3PPF_Conclusion}
summarizes the work with concluded remarks.

\section{Math Notation \label{sec:SO3PPF_Math-Notations}}

In this paper, $\mathbb{R}_{+}$ refers to the set of non-negative
real numbers. $\mathbb{R}^{n}$ is the real space with $n$ dimensions
while $\mathbb{R}^{n\times m}$ stands for the real space of dimensions
$n\times m$. The Euclidean norm of $x\in\mathbb{R}^{n}$ is expressed
as $||x||=\sqrt{x^{\top}x}$, with $^{\top}$ denoting the transpose
of the associated component. $\lambda\left(\cdot\right)$ represents
a group of eigenvalues of a matrix while $\underline{\lambda}\left(\cdot\right)$
is the minimum eigenvalue. $\mathbf{I}_{n}$ denotes an $n$-by-$n$
identity matrix, and zero vector $\underline{\boldsymbol{0}}_{n}$
has $n$-rows and one column. Let $\mathbb{SO}\left(3\right)$ represent
the Special Orthogonal Group. The rigid-body attitude is expressed
as a rotational matrix $R$:
\[
\mathbb{SO}\left(3\right):=\left\{ \left.R\in\mathbb{R}^{3\times3}\right|R^{\top}R=\mathbf{I}_{3}\text{, }{\rm det}\left(R\right)=1\right\} 
\]
where $\mathbf{I}_{n}$ is an $n$-by-$n$ identity matrix, and ${\rm det\left(\cdot\right)}$
denotes the determinant of the matrix. $\mathfrak{so}\left(3\right)$
is the Lie-algebra associated with $\mathbb{SO}\left(3\right)$ and
can be defined by
\[
\mathfrak{so}\left(3\right):=\left\{ \left.\mathcal{A}\in\mathbb{R}^{3\times3}\right|\mathcal{A}^{\top}=-\mathcal{A}\right\} 
\]
where $\mathcal{A}$ is the space of skew-symmetric matrices. Define
the map $\left[\cdot\right]_{\times}:\mathbb{R}^{3}\rightarrow\mathfrak{so}\left(3\right)$
such that
\[
\mathcal{A}=\left[\alpha\right]_{\times}=\left[\begin{array}{ccc}
0 & -\alpha_{3} & \alpha_{2}\\
\alpha_{3} & 0 & -\alpha_{1}\\
-\alpha_{2} & \alpha_{1} & 0
\end{array}\right],\hspace{1em}\alpha=\left[\begin{array}{c}
\alpha_{1}\\
\alpha_{2}\\
\alpha_{3}
\end{array}\right]\in\mathbb{R}^{3}
\]
For all $\alpha,\beta\in\mathbb{R}^{3}$, we have $\left[\alpha\right]_{\times}\beta=\alpha\times\beta$
such that the cross product of two vectors is denoted by $\times$.
Consider that the vex operator is the inverse of $\left[\cdot\right]_{\times}$,
represented by $\mathbf{vex}:\mathfrak{so}\left(3\right)\rightarrow\mathbb{R}^{3}$
where $\mathbf{vex}\left(\mathcal{A}\right)=\alpha$ for all $\alpha\in\mathbb{R}^{3}$
and $\mathcal{A}\in\mathfrak{so}\left(3\right)$. Let $\boldsymbol{\mathcal{P}}_{a}$
stand for the anti-symmetric projection component on the Lie-algebra
$\mathfrak{so}\left(3\right)$ \cite{murray1994mathematical}, expressed
as $\boldsymbol{\mathcal{P}}_{a}:\mathbb{R}^{3\times3}\rightarrow\mathfrak{so}\left(3\right)$,
thus
\[
\boldsymbol{\mathcal{P}}_{a}\left(\mathcal{B}\right)=\frac{1}{2}\left(\mathcal{B}-\mathcal{B}^{\top}\right)\in\mathfrak{so}\left(3\right)
\]
for all $\mathcal{B}\in\mathbb{R}^{3\times3}$. The normalized Euclidean
distance of a rotation matrix on $\mathbb{SO}\left(3\right)$ can
be represented as follows
\begin{equation}
||R||_{I}:=\frac{1}{4}{\rm Tr}\left\{ \mathbf{I}_{3}-R\right\} \label{eq:SO3PPF_Ecul_Dist}
\end{equation}
with ${\rm Tr}\left\{ \cdot\right\} $ being the trace of the associated
matrix and $||R||_{I}\in\left[0,1\right]$. Knowledge of axis parameterization
$u\in\mathbb{R}^{3}$ and angle of rotation $\alpha\in\mathbb{R}$
is sufficient for the reconstruction of the rigid-body attitude. This
attitude reconstruction method is referred to as angle-axis parameterization
\cite{shuster1993survey}. One can define the mapping of angle-axis
parameterization to $\mathbb{SO}\left(3\right)$ by $\mathcal{R}_{\alpha}:\mathbb{R}\times\mathbb{R}^{3}\rightarrow\mathbb{SO}\left(3\right)$
and obtain
\begin{align}
\mathcal{R}_{\alpha}\left(\alpha,u\right) & =\exp\left(-\alpha\left[u\right]_{\times}\right)\nonumber \\
& =\mathbf{I}_{3}+\sin\left(\alpha\right)\left[u\right]_{\times}+\left(1-\cos\left(\alpha\right)\right)\left[u\right]_{\times}^{2}\label{eq:SO3PPF_att_ang}
\end{align}
The identities below will be used in the filter derivation

\begin{align}
\left[\alpha\times\beta\right]_{\times}= & \beta\alpha^{\top}-\alpha\beta^{\top},\quad\alpha,\beta\in{\rm \mathbb{R}}^{3}\label{eq:SO3PPF_Identity1}\\
\left[R\alpha\right]_{\times}= & R\left[\alpha\right]_{\times}R^{\top},\quad R\in\mathbb{SO}\left(3\right),\alpha\in\mathbb{R}^{3}\label{eq:SO3PPF_Identity2}\\
\left[\alpha\right]_{\times}^{2}= & -\alpha^{\top}\alpha\mathbf{I}_{3}+\alpha\alpha^{\top},\quad\alpha\in\mathbb{R}^{3}\label{eq:SO3PPF_Identity3}\\
B\left[\alpha\right]_{\times}+\left[\alpha\right]_{\times}B= & {\rm Tr}\left\{ B\right\} \left[\alpha\right]_{\times}-\left[B\alpha\right]_{\times},\nonumber \\
& \qquad B=B^{\top}\in\mathbb{R}^{3\times3},\alpha\in\mathbb{R}^{3}\label{eq:SO3PPF_Identity4}\\
{\rm Tr}\left\{ \left[A,B\right]\right\} = & {\rm Tr}\left\{ AB-BA\right\} =0,\quad A,B\in\mathbb{R}^{3\times3}\label{eq:SO3PPF_Identity5}\\
{\rm Tr}\left\{ B\left[\alpha\right]_{\times}\right\} = & 0,\quad B=B^{\top}\in\mathbb{R}^{3\times3},\alpha\in\mathbb{R}^{3}\label{eq:SO3PPF_Identity6}
\end{align}
\begin{align}
{\rm Tr}\left\{ A\left[\alpha\right]_{\times}\right\} = & {\rm Tr}\left\{ \boldsymbol{\mathcal{P}}_{a}\left(A\right)\left[\alpha\right]_{\times}\right\} =-2\mathbf{vex}\left(\boldsymbol{\mathcal{P}}_{a}\left(A\right)\right)^{\top}\alpha,\nonumber \\
& \qquad A\in\mathbb{R}^{3\times3},\alpha\in\mathbb{R}^{3}\label{eq:SO3PPF_Identity7}
\end{align}

\section{Problem Formulation with Prescribed Performance \label{sec:SO3PPF_Problem-Formulation-in}}

Attitude estimator relies on a collection of inertial-frame and body-frame
vectorial measurements. In this section,
the attitude problem is defined, and body-frame and gyroscope measurements
are presented. Next, the attitude error is defined and reformulated
to satisfy a desired measure of transient and steady-state performance.

\subsection{Attitude Kinematics and Measurements\label{subsec:SO3PPF_Attitude-Kinematics}}

\begin{figure}[h]
	\centering{}\includegraphics[scale=0.7]{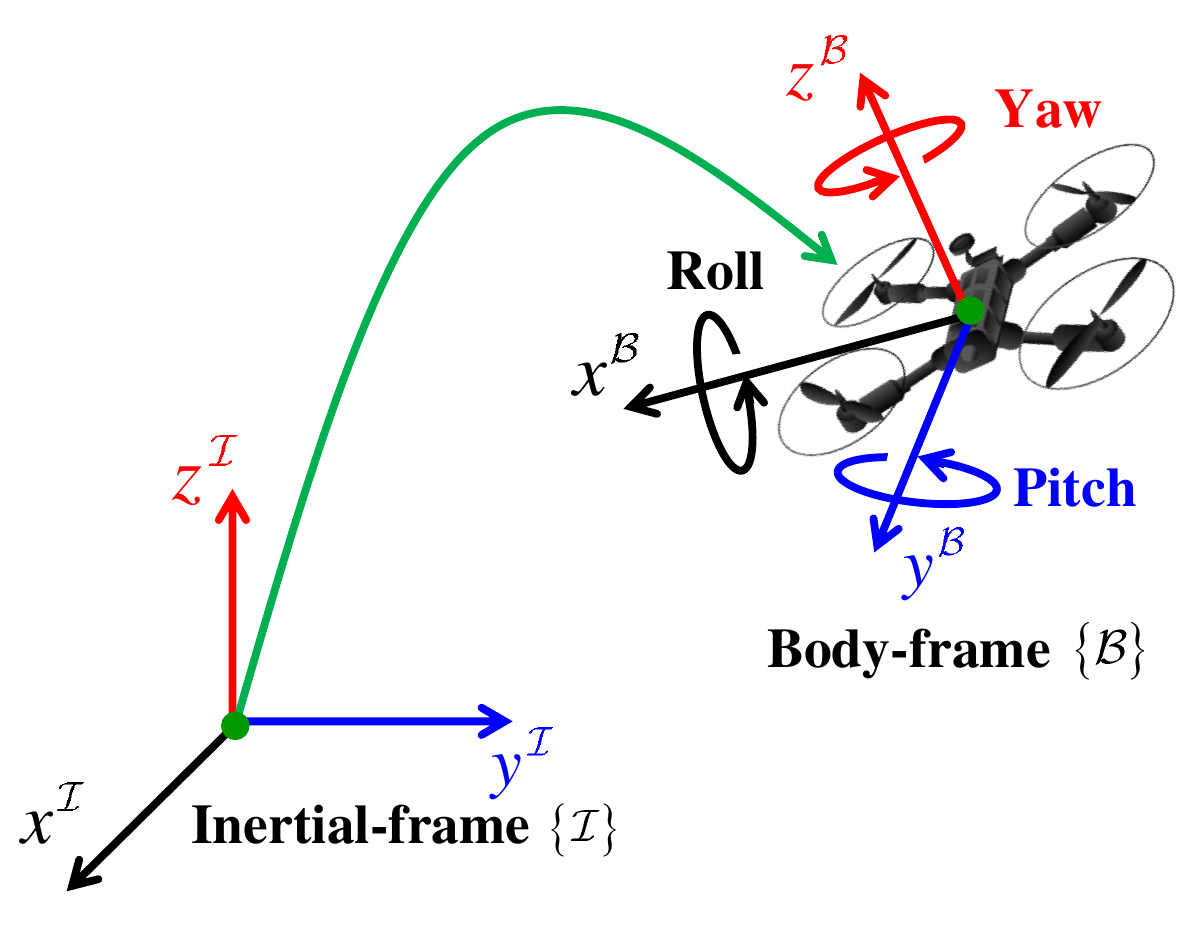}\caption{The relative orientation between body-frame and inertial-frame of
		a rigid-body in 3D space.}
	\label{fig:SO3PPF_1} 
\end{figure}

$R\in\mathbb{SO}\left(3\right)$ stands for the rotational matrix,
and therefore the orientation of the object in the body-frame $\left\{ \mathcal{B}\right\} $
relative to the inertial-frame $\left\{ \mathcal{I}\right\} $ can
be represented by the attitude matrix $R\in\left\{ \mathcal{B}\right\} $
as illustrated in Figure \ref{fig:SO3PPF_1}. 

Let the superscripts $\mathcal{I}$ and $\mathcal{B}$ denote a vector
associated with the inertial-frame and body-frame, respectively. Consider
${\rm v}_{i}^{\mathcal{I}}\in\mathbb{R}^{3}$ to be a known vector
in the inertial-frame and to be measured in the coordinate system
fixed to the rigid-body such that

\begin{equation}
{\rm v}_{i}^{\mathcal{B}}=R^{\top}{\rm v}_{i}^{\mathcal{I}}+{\rm b}_{i}^{\mathcal{B}}+\omega_{i}^{\mathcal{B}}\label{eq:SO3PPF_Vect_True}
\end{equation}
where ${\rm v}_{i}^{\mathcal{B}}\in\mathbb{R}^{3}$ is the $i$th
body-frame measurement associated with ${\rm v}_{i}^{\mathcal{I}}$.
${\rm b}_{i}^{\mathcal{B}}\in\mathbb{R}^{3}$ stands for the bias
component, and $\omega_{i}^{\mathcal{B}}\in\mathbb{R}^{3}$ denotes
the noise component attached to the $i$th body-frame measurement
for $i=1,2,\ldots,n$. Suppose that the instantaneous set of size
$n\geq2$ consisting of known inertial-frame and measured body-frame
vectors is non-collinear. Therefore, the attitude can be established.
Moreover, two non-collinear vectors $\left(n=2\right)$ are generally
sufficient for attitude reconstruction, e.g., \cite{shuster1981three,hashim2018SO3Stochastic,crassidis2007survey,hashim2018Conf1,mahony2008nonlinear,hashim2018SE3Stochastic}.
In case when $n=2$, the third inertial-frame and body-frame vectors
can be obtained by the cross product such that ${\rm v}_{3}^{\mathcal{I}}={\rm v}_{1}^{\mathcal{I}}\times{\rm v}_{2}^{\mathcal{I}}$
and ${\rm v}_{3}^{\mathcal{B}}={\rm v}_{1}^{\mathcal{B}}\times{\rm v}_{2}^{\mathcal{B}}$,
respectively. The inertial-frame and body-frame vectors can be normalized
and their normalized values can be implemented in the estimation of
the attitude in the following manner
\begin{equation}
\upsilon_{i}^{\mathcal{I}}=\frac{{\rm v}_{i}^{\mathcal{I}}}{||{\rm v}_{i}^{\mathcal{I}}||},\hspace{1em}\upsilon_{i}^{\mathcal{B}}=\frac{{\rm v}_{i}^{\mathcal{B}}}{||{\rm v}_{i}^{\mathcal{B}}||}\label{eq:SO3PPF_Vector_norm}
\end{equation}
Hence, the attitude can be obtained knowing $\upsilon_{i}^{\mathcal{I}}$
and $\upsilon_{i}^{\mathcal{B}}$. For simplicity, it is considered
that the body frame vector (${\rm v}_{i}^{\mathcal{B}}$) is noise
and bias free in the stability analysis. The Simulation Section, on
the other hand, takes noise and bias associated with the measurements
into account. The angular velocity of the moving body relative to
the inertial-frame is measured by the rate gyros as
\begin{equation}
\Omega_{m}=\Omega+b+\omega\label{eq:SO3PPF_Angular}
\end{equation}
where $\Omega\in\mathbb{R}^{3}$ is the true value of angular velocity
and $b$ and $\omega$ denote the bias and noise components, respectively,
attached to the measurement of angular velocity for all $b,\omega\in\mathbb{R}^{3}$.
The kinematics of the true attitude are described by
\begin{equation}
\dot{R}=R\left[\Omega\right]_{\times}\label{eq:SO3PPF_R_dynam}
\end{equation}
where $\Omega\in\left\{ \mathcal{B}\right\} $. Considering the normalized
Euclidean distance of $R$ in \eqref{eq:SO3PPF_Ecul_Dist} and the
identity in \eqref{eq:SO3PPF_Identity7}, the kinematics of the true
attitude in \eqref{eq:SO3PPF_R_dynam} can be defined in terms of
normalized Euclidean distance as
\begin{align}
\frac{d}{dt}||R||_{I} & =-\frac{1}{4}{\rm Tr}\left\{ \dot{R}\right\} \nonumber \\
& =-\frac{1}{4}{\rm Tr}\left\{ \boldsymbol{\mathcal{P}}_{a}\left(R\right)\left[\Omega\right]_{\times}\right\} \nonumber \\
& =\frac{1}{2}\mathbf{vex}\left(\boldsymbol{\mathcal{P}}_{a}\left(R\right)\right)^{\top}\Omega\label{eq:SO3PPF_NormR_dynam}
\end{align}
For the sake of simplicity, let us neglect the noise attached to angular
velocity measurements such that the kinematics of the normalized Euclidean
distance in \eqref{eq:SO3PPF_NormR_dynam} become
\begin{equation}
\frac{d}{dt}||R||_{I}=\frac{1}{2}\mathbf{vex}\left(\boldsymbol{\mathcal{P}}_{a}\left(R\right)\right)^{\top}\left(\Omega_{m}-b\right)\label{eq:SO3PPF_NormR_Bias}
\end{equation}
Now, we introduce Lemma \ref{Lemm:SO3PPF_1} which is going to be
applicable in the subsequent filter derivation. 
\begin{lem}
	\label{Lemm:SO3PPF_1}Let $R\in\mathbb{SO}\left(3\right)$, $M^{\mathcal{B}}=\left(M^{\mathcal{B}}\right)^{\top}\in\mathbb{R}^{3\times3}$,
	$M^{\mathcal{B}}$ be nonsingular, ${\rm Tr}\left\{ M^{\mathcal{B}}\right\} =3$,
	and $\bar{\mathbf{M}}^{\mathcal{B}}={\rm Tr}\left\{ M^{\mathcal{B}}\right\} \mathbf{I}_{3}-M^{\mathcal{B}}$,
	while the minimum singular value of $\bar{\mathbf{M}}^{\mathcal{B}}$
	is $\underline{\lambda}:=\underline{\lambda}\left(\bar{\mathbf{M}}^{\mathcal{B}}\right)$.
	Then, the following holds:
	\begin{equation}
	||\mathbf{vex}\left(\boldsymbol{\mathcal{P}}_{a}\left(R\right)\right)||^{2}=4\left(1-||R||_{I}\right)||R||_{I}\label{eq:SO3PPF_lemm1_1}
	\end{equation}
	\begin{align}
	\frac{2}{\underline{\lambda}}\frac{||\mathbf{vex}\left(\boldsymbol{\mathcal{P}}_{a}\left(M^{\mathcal{B}}R\right)\right)||^{2}}{1+{\rm Tr}\left\{ \left(M^{\mathcal{B}}\right)^{-1}M^{\mathcal{B}}R\right\} } & \geq\left\Vert M^{\mathcal{B}}R\right\Vert _{I}\label{eq:SO3PPF_lemm1_3}
	\end{align}
	\textbf{Proof. See \nameref{sec:SO3PPF_AppendixA}.} 
\end{lem}

\subsection{Estimator Structure and Error Criteria}

The goal of attitude estimator in this work is to achieve accurate
estimate of the true attitude satisfying transient as well as steady-state
characteristics. In this subsection general framework of the nonlinear
attitude filter on $\mathbb{SO}\left(3\right)$ is introduced. Next,
the error dynamics are expressed with respect to normalized Euclidean
distance. Let $\hat{R}$ denote the estimate of the true attitude
$R$ and $\tilde{R}=R^{\top}\hat{R}$ denote the attitude error between body-frame and estimator-frame. Consider the following nonlinear attitude filter on $\mathbb{SO}\left(3\right)$
\begin{align}
\dot{\hat{R}} & =\hat{R}\left[\Omega_{m}-\hat{b}-W\right]_{\times},\quad\hat{R}\left(0\right)=\hat{R}_{0}\label{eq:SO3PPF_Rest_dot}\\
\dot{\hat{b}} & =\frac{1}{2}\mathbf{K}_{b}\mathbf{vex}\left(\boldsymbol{\mathcal{P}}_{a}\left(\boldsymbol{\Phi}\right)\right),\quad\hat{b}\left(0\right)=\hat{b}_{0}\label{eq:SO3PPF_b_est}\\
W & =2\mathbf{K}_{W}\mathbf{vex}\left(\boldsymbol{\mathcal{P}}_{a}\left(\boldsymbol{\Phi}\right)\right)\label{eq:SO3PPF_W}
\end{align}
with $\hat{b}$ being the estimate of the true rate-gyro bias $b$,
$\mathbf{K}_{b}$ being a time-variant gain associated with $\hat{b}$,
$\mathbf{K}_{W}$ being a time-variant gain associated with the correction
factor $W$, and $\boldsymbol{\Phi}$ being a matrix associated with
attitude error $\tilde{R}$. Define the unstable set $\mathcal{U}\subseteq\mathbb{SO}\left(3\right)$
by $\mathcal{U}:=\left\{ \left.\tilde{R}_{0}\right|{\rm Tr}\left\{ \tilde{R}_{0}\right\} =-1\right\} $
with $\tilde{R}_{0}=\tilde{R}\left(0\right)$. $\mathbf{K}_{b}$,
$\mathbf{K}_{W}$, and $\boldsymbol{\Phi}$ will be defined subsequently.
In particular, the dynamic gains $\mathbf{K}_{b}$ and $\mathbf{K}_{W}$
will be selected such that their values become increasingly aggressive
as $\tilde{R}$ approaches the unstable equilibria ${\rm Tr}\left\{ \tilde{R}_{0}\right\} \rightarrow-1$,
and reduce significantly as $\tilde{R}$ approaches $\mathbf{I}_{3}$.
\begin{rem}
	\label{rem:RemNew}In the conventional design of nonlinear attitude
	filters, for example \cite{crassidis2007survey,mahony2008nonlinear,hamel2006attitude,grip2012attitude},
	$\mathbf{K}_{b}$ and $\mathbf{K}_{W}$ are selected as positive constant
	gains. However, the weakness of the conventional design of nonlinear
	attitude filters is that smaller values of $\mathbf{K}_{b}$ and $\mathbf{K}_{W}$
	result in slower transient performance with less oscillatory behavior
	in the steady-state. In contrast, higher values of $\mathbf{K}_{b}$
	and $\mathbf{K}_{W}$ generate faster transient performance with higher
	oscillation in the steady-state. 
\end{rem}
Consider the error between body-frame and estimator-frame being defined
as
\begin{equation}
\tilde{R}=R^{\top}\hat{R}\label{eq:SO3PPF_R_error}
\end{equation}
Also, define the error in bias estimation by
\begin{align}
\tilde{b} & =b-\hat{b}\label{eq:SO3PPF_b_tilde}
\end{align}

\noindent From \eqref{eq:SO3PPF_R_dynam} and \eqref{eq:SO3PPF_Rest_dot}
the error dynamics can be found to be
\begin{align}
\dot{\tilde{R}} & =R^{\top}\hat{R}\left[\Omega_{m}-\hat{b}-W\right]_{\times}-\left[\Omega\right]_{\times}R^{\top}\hat{R}\nonumber \\
& =\tilde{R}\left[\tilde{b}-W\right]_{\times}+\tilde{R}\left[\Omega\right]_{\times}-\left[\Omega\right]_{\times}\tilde{R}\label{eq:SO3PPF_Rtilde_dot}
\end{align}
Considering \eqref{eq:SO3PPF_R_dynam} and \eqref{eq:SO3PPF_NormR_dynam},
the error dynamics in \eqref{eq:SO3PPF_Rtilde_dot} are represented
with regards to normalized Euclidean distance 
\begin{align}
\frac{d}{dt}||\tilde{R}||_{I} & =\frac{d}{dt}\frac{1}{4}{\rm Tr}\left\{ \mathbf{I}_{3}-\tilde{R}\right\} \nonumber \\
& =-\frac{1}{4}{\rm Tr}\left\{ \tilde{R}\left[\tilde{b}-W\right]_{\times}\right\} -\frac{1}{4}{\rm Tr}\left\{ \left[\tilde{R},\left[\Omega\right]_{\times}\right]\right\} \nonumber \\
& =\frac{1}{2}\mathbf{vex}\left(\boldsymbol{\mathcal{P}}_{a}\left(\tilde{R}\right)\right)^{\top}\left(\tilde{b}-W\right)\label{eq:SO3PPF_NormRtilde_dot}
\end{align}
where ${\rm Tr}\left\{ \tilde{R}\left[\tilde{b}\right]_{\times}\right\} =-2\mathbf{vex}\left(\boldsymbol{\mathcal{P}}_{a}\left(\tilde{R}\right)\right)^{\top}\tilde{b}$
as given in \eqref{eq:SO3PPF_Identity7} and ${\rm Tr}\left\{ \left[\tilde{R},\left[\Omega\right]_{\times}\right]\right\} =0$
as defined in \eqref{eq:SO3PPF_Identity5}.

\subsection{Prescribed Performance \label{subsec:SO3PPF_Prescribed-Performance}}

\begin{figure*}[!h]
	\centering{}\includegraphics[scale=0.42]{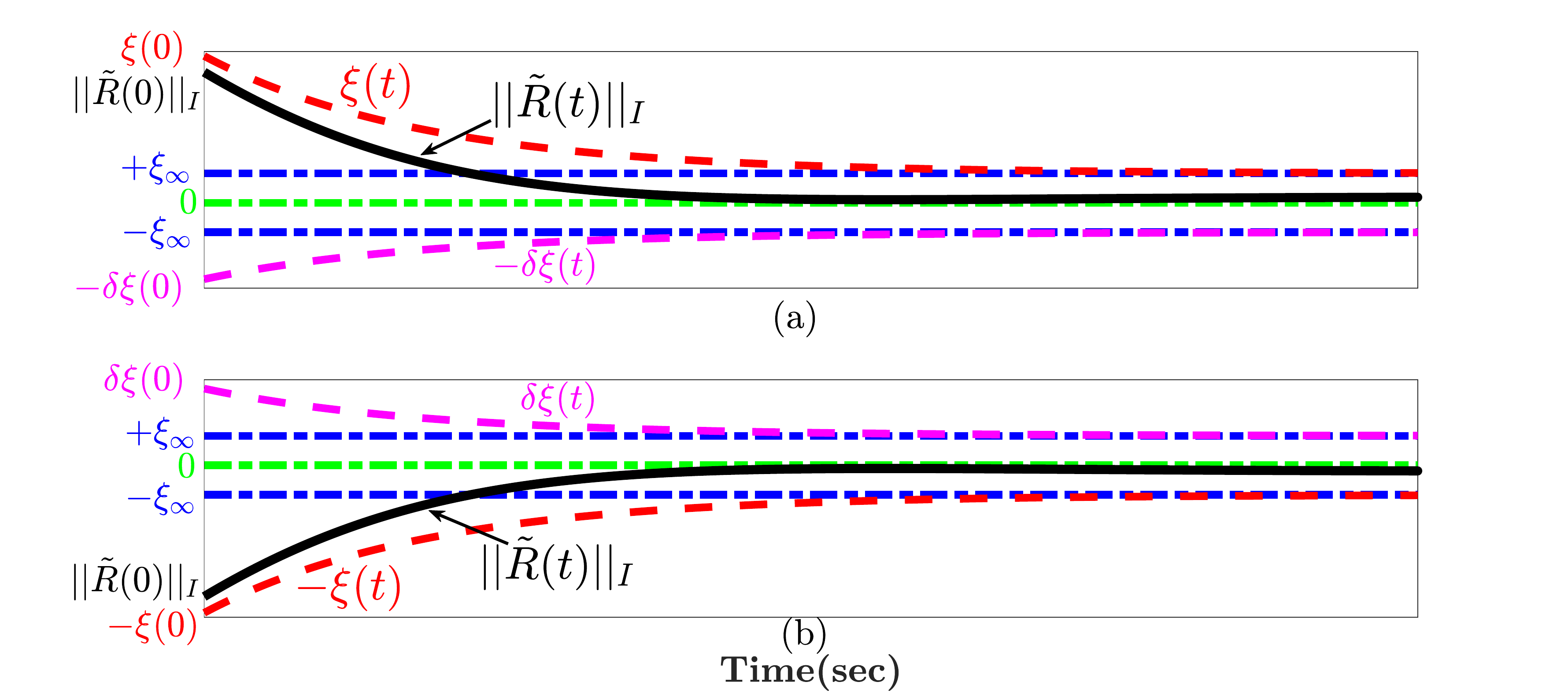} \caption{A detailed representation of tracking normalized Euclidean distance
		error $||\tilde{R}\left(t\right)||_{I}$ with PPF satisfying (a) Eq.
		\eqref{eq:SO3PPF_ePos}; (b) Eq. \eqref{eq:SO3PPF_eNeg}.}
	\label{fig:SO3PPF_2} 
\end{figure*}

This subsection aims to reformulate the problem such that the normalized
Euclidean distance of the attitude error $||\tilde{R}\left(t\right)||_{I}$
satisfies the predefined transient as well as steady-state measures
set by the user. Initially, the error $||\tilde{R}\left(t\right)||_{I}$
is contained within a predefined large set and decreases systematically
and smoothly to a predefined narrow set through a prescribed performance
function (PPF) \cite{bechlioulis2008robust}. This is accomplished
by first defining a configuration error function \cite{bechlioulis2008robust,hashim2017neuro,hashim2017adaptive}.
Let $\xi\left(t\right)$ be a positive smooth and time-decreasing
performance function such that $\xi:\mathbb{R}_{+}\to\mathbb{R}_{+}$
and $\lim\limits _{t\to\infty}\xi\left(t\right)=\xi_{\infty}>0$.
The general expression of the PPF is as follows
\begin{equation}
\xi\left(t\right)=\left(\xi_{0}-\xi_{\infty}\right)\exp\left(-\ell t\right)+\xi_{\infty}\label{eq:SO3PPF_Presc}
\end{equation}
where $\xi_{0}=\xi\left(0\right)$ is the upper bound of the predefined
large set, also known to be the initial value of the PPF, $\xi_{\infty}$
is the upper bound of the small set such that the steady-state error
is confined by $\pm\xi_{\infty}$, while $\ell$ is a positive constant
controlling the convergence rate of the set boundaries $\xi\left(t\right)$
with respect to time from $\xi_{0}$ to $\xi_{\infty}$. It is sufficient
to force $||\tilde{R}\left(t\right)||_{I}$ to obey a predefined transient
and steady-state characteristics, if the following conditions are
met:
\begin{align}
-\delta\xi\left(t\right)<||\tilde{R}\left(t\right)||_{I}<\xi\left(t\right), & \text{ if }||\tilde{R}\left(0\right)||_{I}\geq0,\forall t\geq0\label{eq:SO3PPF_ePos}\\
-\xi\left(t\right)<||\tilde{R}\left(t\right)||_{I}<\delta\xi\left(t\right), & \text{ if }||\tilde{R}\left(0\right)||_{I}<0,\forall t\geq0\label{eq:SO3PPF_eNeg}
\end{align}
where $\delta$ is selected such that $1\geq\delta\geq0$. The tracking
error $||\tilde{R}\left(t\right)||_{I}$, with PPF decreasing systematically
from a known large set to a known small set in accordance with \eqref{eq:SO3PPF_ePos}
and \eqref{eq:SO3PPF_eNeg} is illustrated in Figure \ref{fig:SO3PPF_2}.

\begin{rem}
	\label{rem3} As explained in \cite{bechlioulis2008robust,hashim2017neuro},
	knowing the sign of $||\tilde{R}\left(0\right)||_{I}$ is sufficient
	to satisfy the performance constraints and maintain the error convergence
	within the predefined dynamically decreasing boundaries for all $t>0$.
	Since $||\tilde{R}\left(0\right)||_{I}\in\left[0,1\right]$, $||\tilde{R}\left(0\right)||_{I}$
	is guaranteed to be greater than or equal to 0 for any attitude initialization,
	and therefore the only possible condition is \eqref{eq:SO3PPF_ePos}.
	If the condition in \eqref{eq:SO3PPF_ePos} is met, the maximum steady-state
	error will be less than $\xi_{\infty}$, the maximum overshoot will
	be less than $-\delta\xi\left(0\right)$, and $||\tilde{R}\left(t\right)||_{I}$
	will be confined between $\xi\left(t\right)$ and $\delta\xi\left(t\right)$
	as given in the upper portion of Figure \ref{fig:SO3PPF_2}. 
\end{rem}
Let us define
\begin{equation}
||\tilde{R}\left(t\right)||_{I}=\xi\left(t\right)\mathcal{Z}\left(\mathcal{E}\right)\label{eq:SO3PPF_e_Trans}
\end{equation}
with $\xi\left(t\right)\in\mathbb{R}$ being given in \eqref{eq:SO3PPF_Presc},
$\mathcal{E}\in\mathbb{R}$ being a transformed error, and $\mathcal{Z}\left(\mathcal{E}\right)$
being a smooth function which satisfies Assumption \ref{Assum:SO3PPF_1}: 
\begin{assum}
	\label{Assum:SO3PPF_1}The smooth function $\mathcal{Z}\left(\mathcal{E}\right)$
	must satisfy \cite{bechlioulis2008robust}:
\begin{enumerate}
	\item[P 1)] $\mathcal{Z}\left(\mathcal{E}\right)$ is smooth and strictly increasing. 
	\item[P 2)] $\mathcal{Z}\left(\mathcal{E}\right)$ is bounded between two predefined
	bounds \\
	$-\underline{\delta}<\mathcal{Z}\left(\mathcal{E}\right)<\bar{\delta},{\rm \text{for}}||\tilde{R}\left(0\right)||_{I}\geq0$\\
	with $\bar{\delta}$ and $\underline{\delta}$ being positive constants
	and $\underline{\delta}\leq\bar{\delta}$. 
	\item[P 3)] $\underset{\mathcal{E}\rightarrow-\infty}{\lim}\mathcal{Z}\left(\mathcal{E}\right)=-\underline{\delta}$
	and $\underset{\mathcal{E}\rightarrow+\infty}{\lim}\mathcal{Z}\left(\mathcal{E}\right)=\bar{\delta}$
	where 
	\begin{equation}
	\mathcal{Z}\left(\mathcal{E}\right)=\frac{\bar{\delta}\exp\left(\mathcal{E}\right)-\underline{\delta}\exp\left(-\mathcal{E}\right)}{\exp\left(\mathcal{E}\right)+\exp\left(-\mathcal{E}\right)}\label{eq:SO3PPF_Smooth}
	\end{equation}
\end{enumerate}
\end{assum}
One could find the transformed error to be
\begin{equation}
\mathcal{E}\left(||\tilde{R}\left(t\right)||_{I},\xi\left(t\right)\right)=\mathcal{Z}^{-1}\left(\frac{||\tilde{R}\left(t\right)||_{I}}{\xi\left(t\right)}\right)\label{eq:SO3PPF_Trans1}
\end{equation}
where $\mathcal{E}\in\mathbb{R}$, $\mathcal{Z}\in\mathbb{R}$ and
$\mathcal{Z}^{-1}\in\mathbb{R}$ are smooth functions. For clarity,
let $\xi:=\xi\left(t\right)$, $||\tilde{R}||_{I}:=||\tilde{R}\left(t\right)||_{I}$
and $\mathcal{E}:=\mathcal{E}\left(\cdot,\cdot\right)$. The transformed
error $\mathcal{E}$ plays a prominent role driving the error dynamics
from constrained form in either \eqref{eq:SO3PPF_ePos} or \eqref{eq:SO3PPF_eNeg}
to that in \eqref{eq:SO3PPF_Trans1} which is unconstrained. One can
find from \eqref{eq:SO3PPF_Smooth} that the transformed error is
\begin{equation}
\begin{aligned}\mathcal{E}= & \frac{1}{2}\text{ln}\frac{\underline{\delta}+||\tilde{R}||_{I}/\xi}{\bar{\delta}-||\tilde{R}||_{I}/\xi}\end{aligned}
\label{eq:SO3PPF_trans3}
\end{equation}

\begin{rem}
	\label{rem:SO3PPF_1} \cite{bechlioulis2008robust,hashim2017neuro}
	Consider the transformed error in \eqref{eq:SO3PPF_trans3}. If $\mathcal{E}\left(t\right)$
	is guaranteed to be bounded for all $t\geq0$, the performance function
	$\xi\left(t\right)$ can be used to bound the transient and steady-state
	of the tracking error ($||\tilde{R}||_{I}$) allowing it to achieve
	the prescribed performance.
\end{rem}
\begin{prop}
	\label{Prop:SO3PPF_1}Consider the normalized Euclidean distance error
	$||\tilde{R}||_{I}$ being defined by \eqref{eq:SO3PPF_Ecul_Dist}
	and from \eqref{eq:SO3PPF_e_Trans}, \eqref{eq:SO3PPF_Smooth}, \eqref{eq:SO3PPF_Trans1}
	let the transformed error be given as in \eqref{eq:SO3PPF_trans3}
	with $\underline{\delta}=\bar{\delta}$. Then the following statements
	hold.
\begin{enumerate}
	\item[(i)] The transformed error $\mathcal{E}>0\forall||\tilde{R}||_{I}\neq0$
	and $\mathcal{E}=0$ only at $||\tilde{R}||_{I}=0$. 
	\item[(ii)] The critical point of $\mathcal{E}$ satisfies $||\tilde{R}||_{I}=0$. 
	\item[(iii)] The only critical point of $\mathcal{E}$ is $\tilde{R}=\mathbf{I}_{3}$. 
\end{enumerate}
\end{prop}
\textbf{Proof.} Letting $\underline{\delta}=\bar{\delta}$ with the
prescribed performance constraints $||\tilde{R}||_{I}\leq\xi$, the
expression $\left(\underline{\delta}+||\tilde{R}||_{I}/\xi\right)/\left(\bar{\delta}-||\tilde{R}||_{I}/\xi\right)$
in \eqref{eq:SO3PPF_trans3} is always greater than or equal to 1.
Accordingly, $\mathcal{E}>0\forall||\tilde{R}||_{I}\neq0$ and $\mathcal{E}=0$
at $||\tilde{R}||_{I}=0$ which proves (i). For (ii) and (iii), from
\eqref{eq:SO3PPF_Ecul_Dist}, $||\tilde{R}||_{I}=0$ if and only if
$\tilde{R}=\mathbf{I}_{3}$. Thus, the critical point of $\mathcal{E}$
satisfies $\tilde{R}=\mathbf{I}_{3}$ and, consequently, also satisfies
$||\tilde{R}||_{I}=0$ which proves (ii) and (iii). Let us define
a new variable $\mu:=\mu\left(||\tilde{R}||_{I},\xi\right)$ such
that
\begin{equation}
\begin{split}\mu & =\frac{1}{2\xi}\frac{\partial\mathcal{Z}^{-1}\left(||\tilde{R}||_{I}/\xi\right)}{\partial\left(||\tilde{R}||_{I}/\xi\right)}\\
& =\frac{1}{2\xi}\left(\frac{1}{\underline{\delta}+||\tilde{R}||_{I}/\xi}+\frac{1}{\bar{\delta}-||\tilde{R}||_{I}/\xi}\right)
\end{split}
\label{eq:SO3PPF_Trans_Aux}
\end{equation}
Consequently, the derivative of the transformed error is governed
b{\small{}y
	\begin{align}
	\dot{\mathcal{E}} & =\frac{1}{2\xi}\left(\frac{1}{\underline{\delta}+||\tilde{R}||_{I}/\xi}+\frac{1}{\bar{\delta}-||\tilde{R}||_{I}/\xi}\right)\left(\frac{d}{dt}||\tilde{R}||_{I}-\frac{\dot{\xi}}{\xi}||\tilde{R}||_{I}\right)\nonumber \\
	& =\mu\left(\frac{1}{2}\mathbf{vex}\left(\boldsymbol{\mathcal{P}}_{a}\left(\tilde{R}\right)\right)^{\top}\left(\tilde{b}-W\right)-\frac{\dot{\xi}}{\xi}||\tilde{R}||_{I}\right)\label{eq:SO3PPF_Trans_dot}
	\end{align}
}with direct substitution of \eqref{eq:SO3PPF_NormRtilde_dot} in
\eqref{eq:SO3PPF_Trans_dot}. Next section presents two nonlinear
attitude filters on $\mathbb{SO}\left(3\right)$ with prescribed performance
which guarantees $\mathcal{E}\in\mathcal{L}_{\infty},\forall t\geq0$
and, thus, satisfies \eqref{eq:SO3PPF_ePos} provided that $0\leq||\tilde{R}\left(0\right)||_{I}<\xi\left(0\right)$.

\section{Nonlinear Complementary Filters On $\mathbb{SO}\left(3\right)$ with
	Prescribed Performance \label{sec:SO3PPF-Filters}}

The primary objective of this section is to propose two nonlinear
attitude estimators on $\mathbb{SO}\left(3\right)$ with normalized
Euclidean distance error satisfying a predefined transient as well
as steady-state performance given by the user. The constrained error
$||\tilde{R}||_{I}$ is relaxed to unconstrained $\mathcal{E}$. The
first filter is termed a semi-direct attitude filter with prescribed
performance because it requires the attitude to be reconstructed via
the set of vectorial measurements as defined in \eqref{eq:SO3PPF_Vector_norm},
in addition to the measurement of the angular velocity in \eqref{eq:SO3PPF_Angular}.
The second filter is called a direct attitude filter with prescribed
performance because it uses the vectorial measurements in \eqref{eq:SO3PPF_Vector_norm}
and the angular velocity measurement in \eqref{eq:SO3PPF_Angular}
directly without the need for attitude reconstruction.

\subsection{Semi-direct Attitude Filter with Prescribed Performance \label{subsec:SO3PPF_Passive-Filter}}

Let $R_{y}$ denote the reconstructed attitude of $R$. There are
many methods to reconstruct $R_{y}$, for instance, TRIAD \cite{black1964passive},
QUEST \cite{shuster1981three}, or SVD \cite{markley1988attitude}.
Consider the following filter kinematics
\begin{align}
\dot{\hat{R}} & =\hat{R}\left[\Omega_{m}-\hat{b}-W\right]_{\times},\quad\hat{R}\left(0\right)=\hat{R}_{0}\label{eq:SO3PPF_Rest_dot_Ry}\\
\dot{\hat{b}} & =\frac{1}{2}\gamma\mu\mathcal{E}\mathbf{vex}\left(\boldsymbol{\mathcal{P}}_{a}\left(\tilde{R}\right)\right),\quad\hat{b}\left(0\right)=\hat{b}_{0},\tilde{R}=R_{y}^{\top}\hat{R}\label{eq:SO3PPF_b_est_Ry}\\
W & =2\frac{k_{w}\mu\mathcal{E}-\dot{\xi}/4\xi}{1-||\tilde{R}||_{I}}\mathbf{vex}\left(\boldsymbol{\mathcal{P}}_{a}\left(\tilde{R}\right)\right),\hspace{1em}\tilde{R}=R_{y}^{\top}\hat{R}\label{eq:SO3PPF_W_Ry}
\end{align}
with $\mathcal{E}$ and $\mu$ being defined in \eqref{eq:SO3PPF_trans3}
and \eqref{eq:SO3PPF_Trans_Aux}, respectively, $k_{w}$ and $\gamma$
being positive constants, $||\tilde{R}||_{I}=\frac{1}{4}{\rm Tr}\left\{ \mathbf{I}_{3}-\tilde{R}\right\} $
being defined in \eqref{eq:SO3PPF_Ecul_Dist}, $\xi$ being PPF defined
in \eqref{eq:SO3PPF_Presc}, and $\hat{b}$ being the estimate of
$b$. 
\begin{thm}
	\textbf{\label{thm:SO3PPF_1} }Consider the rotation kinematics in
	\eqref{eq:SO3PPF_R_dynam}, measurements of angular velocity in \eqref{eq:SO3PPF_Angular}
	with no noise associated with the measurement $\Omega_{m}=\Omega+b$,
	in addition to other vector measurements given in \eqref{eq:SO3PPF_Vect_True}
	coupled with the filter in \eqref{eq:SO3PPF_Rest_dot_Ry}, \eqref{eq:SO3PPF_b_est_Ry}
	and \eqref{eq:SO3PPF_W_Ry}. Suppose that measurements can be made
	on two or more body-frame non-collinear vectors. Define $\mathcal{U}\subseteq\mathbb{SO}\left(3\right)\times\mathbb{R}^{3}$
	by $\mathcal{U}:=\left\{ \left.\left(\tilde{R}_{0},\tilde{b}_{0}\right)\right|{\rm Tr}\left\{ \tilde{R}_{0}\right\} =-1,\tilde{b}_{0}=\underline{\boldsymbol{0}}_{3}\right\} $
	with $\tilde{R}_{0}=\tilde{R}\left(0\right)$ and $\tilde{b}_{0}=\tilde{b}\left(0\right)$.
	For almost any initial condition such that $\tilde{R}_{0}\notin\mathcal{U}$
	and $\mathcal{E}\left(0\right)\in\mathcal{L}_{\infty}$, then, all
	signals in the closed loop are bounded, $\lim_{t\rightarrow\infty}\mathcal{E}\left(t\right)=0$
	and $\tilde{R}$ asymptotically approaches $\mathbf{I}_{3}$.
\end{thm}
Theorem \ref{thm:SO3PPF_1} guarantees that the observer dynamics
in \eqref{eq:SO3PPF_Rest_dot_Ry}, \eqref{eq:SO3PPF_b_est_Ry} and
\eqref{eq:SO3PPF_W_Ry} are stable with $\mathcal{E}\left(t\right)$
approaching asymptotically the origin. Since, $\mathcal{E}\left(t\right)$
is bounded, $||\tilde{R}||_{I}$ obeys the prescribed transient and
steady-state performance introduced in \eqref{eq:SO3PPF_Presc}.

\textbf{Proof. }Let the error in attitude and bias be defined by $\tilde{R}=R^{\top}\hat{R}$
and $\tilde{b}=b-\hat{b}$ similar to \eqref{eq:SO3PPF_R_error} and
\eqref{eq:SO3PPF_b_tilde}, respectively. From \eqref{eq:SO3PPF_R_dynam}
and \eqref{eq:SO3PPF_Rest_dot_Ry}, the error dynamics can be obtained
as in \eqref{eq:SO3PPF_Rtilde_dot}. Also, in view of \eqref{eq:SO3PPF_R_dynam}
and \eqref{eq:SO3PPF_NormR_dynam}, the error dynamics are analogous
to \eqref{eq:SO3PPF_NormRtilde_dot} in terms of normalized Euclidean
distance. Therefore, considering \eqref{eq:SO3PPF_NormR_dynam} and
\eqref{eq:SO3PPF_Trans_dot}, the derivative of the transformed error
can be found to be
\begin{align}
\dot{\mathcal{E}}= & \mu\left(\frac{1}{2}\mathbf{vex}\left(\boldsymbol{\mathcal{P}}_{a}\left(\tilde{R}\right)\right)^{\top}\left(\tilde{b}-W\right)-\frac{\dot{\xi}}{\xi}||\tilde{R}||_{I}\right)\label{eq:SO3PPF_Trans_dot_Ry}
\end{align}
Consider the following candidate Lyapunov function
\begin{align}
V\left(\mathcal{E},\tilde{b}\right) & =\frac{1}{2}\mathcal{E}^{2}+\frac{1}{2\gamma}||\tilde{b}||^{2}\label{eq:SO3PPF_V_Ry}
\end{align}
Differentiating $V$ in \eqref{eq:SO3PPF_V_Ry} and substituting for
$\dot{\hat{b}}$ and $W$ in \eqref{eq:SO3PPF_b_est_Ry}, and \eqref{eq:SO3PPF_W_Ry},
respectively, one obtains
\begin{align}
\dot{V}= & \mathcal{E}\dot{\mathcal{E}}-\frac{1}{\gamma}\tilde{b}^{\top}\dot{\hat{b}}\nonumber \\
= & \mu\mathcal{E}\left(\frac{1}{2}\mathbf{vex}\left(\boldsymbol{\mathcal{P}}_{a}\left(\tilde{R}\right)\right)^{\top}\left(\tilde{b}-W\right)-\frac{\dot{\xi}}{\xi}||\tilde{R}||_{I}\right)-\frac{1}{\gamma}\tilde{b}^{\top}\dot{\hat{b}}\nonumber \\
= & -\mathcal{E}\mu\left(\frac{k_{w}\mu\mathcal{E}-\dot{\xi}/4\xi}{1-||\tilde{R}||_{I}}\left\Vert \mathbf{vex}\left(\boldsymbol{\mathcal{P}}_{a}\left(\tilde{R}\right)\right)\right\Vert ^{2}+\frac{\dot{\xi}}{\xi}||\tilde{R}||_{I}\right)\label{eq:SO3PPF_Vdot_Ry}
\end{align}
Substituting for $\left\Vert \mathbf{vex}\left(\boldsymbol{\mathcal{P}}_{a}\left(\tilde{R}\right)\right)\right\Vert ^{2}=4\left(1-||\tilde{R}||_{I}\right)||\tilde{R}||_{I}$
as defined in \eqref{eq:SO3PPF_lemm1_1}, the expression in \eqref{eq:SO3PPF_Vdot_Ry}
becomes
\begin{align}
\dot{V}= & -4k_{w}||\tilde{R}||_{I}\mu^{2}\mathcal{E}^{2}\label{eq:SO3PPF_Vdot_Ry_Final}
\end{align}
This implies that $V\left(t\right)\leq V\left(0\right),\forall t\geq0$.
Given $\tilde{R}_{0}\notin\mathcal{U}$ implies that $\tilde{b}$
remains bounded for all $t\geq0$, and, therefore, $\mathcal{E}$
is bounded and well defined for all $t\geq0$. It can be shown that
\begin{align}
\ddot{V}= & -4k_{w}\left(2\left(\mathcal{E}\dot{\mathcal{E}}\mu^{2}+\mathcal{E}^{2}\mu\dot{\mu}\right)||\tilde{R}||_{I}+\mathcal{E}^{2}\mu^{2}||\dot{\tilde{R}}||_{I}\right)\label{eq:SO3PPF_Vddot_Ry}
\end{align}
From \eqref{eq:SO3PPF_Trans_Aux}, it can be found that
\begin{equation}
\dot{\mu}=-\frac{1}{2}\frac{\underline{\delta}\dot{\xi}+||\dot{\tilde{R}}||_{I}}{\left(\underline{\delta}\xi+||\tilde{R}||_{I}\right)^{2}}-\frac{1}{2}\frac{\bar{\delta}\dot{\xi}-||\dot{\tilde{R}}||_{I}}{\left(\bar{\delta}\xi-||\tilde{R}||_{I}\right)^{2}}\label{eq:SO3PPF_mu_dot}
\end{equation}
where $\dot{\xi}=-\ell\left(\xi^{0}-\xi^{\infty}\right)\exp\left(-\ell t\right)$.
Since $||\dot{\tilde{R}}||_{I}$ is bounded, $\dot{\mu}$ is bounded
which shows that $\ddot{V}$ is bounded for all $t\geq0$. Consequently,
$\dot{V}$ is uniformly continuous, and according to Barbalat Lemma,
$\dot{V}\rightarrow0$ indicates that one or more of the following
conditions are true
\begin{enumerate}
	\item $||\mathcal{E}||\rightarrow0$.
	\item $||\tilde{R}||_{I}\rightarrow0$.
	\item $||\mathcal{E}||\rightarrow0$ and $||\tilde{R}||_{I}\rightarrow0$.
\end{enumerate}
as $t\rightarrow\infty$. According to property (i) and (ii) of Proposition
\ref{Prop:SO3PPF_1}, $||\mathcal{E}||\rightarrow0$ means $||\tilde{R}||_{I}\rightarrow0$
and vice versa. Therefore, $\dot{V}\rightarrow0$ as $t\rightarrow\infty$
strictly indicates that $||\mathcal{E}||\rightarrow0$ and $||\tilde{R}||_{I}\rightarrow0$.
As stated by property (iii) of Proposition \ref{Prop:SO3PPF_1}, $||\mathcal{E}||\rightarrow0$
implies that $\tilde{R}$ asymptotically approaches $\mathbf{I}_{3}$.
Hence, $\dot{V}\rightarrow0$ means that $\tilde{R}$ asymptotically
approaches $\mathbf{I}_{3}$, which completes the proof.

\subsection{Direct Attitude Filter with Prescribed Performance\label{subsec:SO3PPF_Explicit-Filter}}

Let $R_{y}$ denote the reconstructed attitude of $R$ obtained through
a set of vectorial measurements in \eqref{eq:SO3PPF_Vector_norm}.
Although there are many methods to reconstruct $R_{y}$, this may
add computational cost. The filter proposed in the previous Subsection
\ref{subsec:SO3PPF_Passive-Filter} requires $R_{y}$ to obtain the
attitude error $\tilde{R}=R_{y}^{\top}\hat{R}$, for example (the
Appendix in \cite{hashim2018SO3Stochastic,hashim2018SE3Stochastic}).
In this Subsection the aforementioned weakness is eliminated by proposing
a nonlinear filter with prescribed performance in terms of direct
measurements from the inertial and body-frame units. Let us recall
$\upsilon_{i}^{\mathcal{I}}\in\left\{ \mathcal{I}\right\} $ and $\upsilon_{i}^{\mathcal{B}}\in\left\{ \mathcal{B}\right\} $
from \eqref{eq:SO3PPF_Vect_True} and \eqref{eq:SO3PPF_Vector_norm}
for $i=1,\ldots,n$. Let us define
\begin{align}
M^{\mathcal{I}} & =\left(M^{\mathcal{I}}\right)^{\top}=\sum_{i=1}^{n}s_{i}\upsilon_{i}^{\mathcal{I}}\left(\upsilon_{i}^{\mathcal{I}}\right)^{\top}\nonumber \\
M^{\mathcal{B}} & =\left(M^{\mathcal{B}}\right)^{\top}=\sum_{i=1}^{n}s_{i}\upsilon_{i}^{\mathcal{B}}\left(\upsilon_{i}^{\mathcal{B}}\right)^{\top}\nonumber \\
& =R^{\top}M^{\mathcal{I}}R\label{eq:SO3PPF_MB_MI}
\end{align}
where $s_{i}>0$ refers to confidence level of the $i$th sensor measurements,
and in this work $s_{i}$ is selected such that $\sum_{i=1}^{n}s_{i}=3$.
According to \eqref{eq:SO3PPF_MB_MI}, $M^{\mathcal{I}}$ and $M^{\mathcal{B}}$
are symmetric matrices. Assume that at least two non-collinear inertial-frame
and measured body-frame vectors are available. If two typical vectors
are available for measurements, $n=2$, the third vector is obtained
by the cross product as mentioned in Subsection \ref{subsec:SO3PPF_Attitude-Kinematics}.
Thereby, the set of vectors is non-collinear and $M^{\mathcal{B}}$
is nonsingular with ${\rm rank}\left(M^{\mathcal{B}}\right)=3$. Hence,
the three eigenvalues of $M^{\mathcal{B}}$ are greater than zero.
Let $\bar{\mathbf{M}}^{\mathcal{B}}={\rm Tr}\left\{ M^{\mathcal{B}}\right\} \mathbf{I}_{3}-M^{\mathcal{B}}\in\mathbb{R}^{3\times3}$,
provided that ${\rm rank}\left(M^{\mathcal{B}}\right)=3$, then, the
following three statements hold (\cite{bullo2004geometric} page.
553): 
\begin{enumerate}
	\item $\bar{\mathbf{M}}^{\mathcal{B}}$ is a symmetric and positive-definite
	matrix. 
	\item The eigenvectors of $M^{\mathcal{B}}$ coincide with the eigenvectors
	of $\bar{\mathbf{M}}^{\mathcal{B}}$. 
	\item Define the three eigenvalues of $M^{\mathcal{B}}$ by $\lambda\left(M^{\mathcal{B}}\right)=\left\{ \lambda_{1},\lambda_{2},\lambda_{3}\right\} $,
	then $\lambda\left(\bar{\mathbf{M}}^{\mathcal{B}}\right)=\{\lambda_{3}+\lambda_{2},\lambda_{3}+\lambda_{1},\lambda_{2}+\lambda_{1}\}$
	such that the minimum singular value $\underline{\lambda}\left(\bar{\mathbf{M}}^{\mathcal{B}}\right)>0$. 
\end{enumerate}
In the remainder of this section, we assume that ${\rm rank}\left(M^{\mathcal{B}}\right)=3$,
and accordingly the three above-mentioned statements are true. Define
\begin{equation}
\hat{\upsilon}_{i}^{\mathcal{B}}=\hat{R}^{\top}\upsilon_{i}^{\mathcal{I}}\label{eq:SO3PPF_vB_hat}
\end{equation}
Define the error in attitude and bias by $\tilde{R}=R^{\top}\hat{R}$
and $\tilde{b}=b-\hat{b}$ which is similar to \eqref{eq:SO3PPF_R_error}
and \eqref{eq:SO3PPF_b_tilde}, respectively. In order to derive the
explicit filter, it is necessary to present the following equations
expressed in terms of vector measurements. From identity \eqref{eq:SO3PPF_Identity1},
one can find

\begin{align*}
\left[\sum_{i=1}^{n}\frac{s_{i}}{2}\hat{\upsilon}_{i}^{\mathcal{B}}\times\upsilon_{i}^{\mathcal{B}}\right]_{\times} & =\sum_{i=1}^{n}\frac{s_{i}}{2}\left(\upsilon_{i}^{\mathcal{B}}\left(\hat{\upsilon}_{i}^{\mathcal{B}}\right)^{\top}-\hat{\upsilon}_{i}^{\mathcal{B}}\left(\upsilon_{i}^{\mathcal{B}}\right)^{\top}\right)\\
& =\frac{1}{2}R^{\top}M^{\mathcal{I}}R\tilde{R}-\frac{1}{2}\tilde{R}^{\top}R^{\top}M^{\mathcal{I}}R\\
& =\boldsymbol{\mathcal{P}}_{a}\left(M^{\mathcal{B}}\tilde{R}\right)
\end{align*}
such that
\begin{equation}
\mathbf{vex}\left(\boldsymbol{\mathcal{P}}_{a}\left(M^{\mathcal{B}}\tilde{R}\right)\right)=\sum_{i=1}^{n}\frac{s_{i}}{2}\hat{\upsilon}_{i}^{\mathcal{B}}\times\upsilon_{i}^{\mathcal{B}}\label{eq:SO3PPF_VEX_VM}
\end{equation}
The normalized Euclidean distance of $M^{\mathcal{B}}\tilde{R}$ can
be found to be
\begin{align}
||M^{\mathcal{B}}\tilde{R}||_{I} & =\frac{1}{4}{\rm Tr}\left\{ \mathbf{I}_{3}-M^{\mathcal{B}}\tilde{R}\right\} \nonumber \\
& =\frac{1}{4}{\rm Tr}\left\{ \mathbf{I}_{3}-\sum_{i=1}^{n}s_{i}\upsilon_{i}^{\mathcal{B}}\left(\hat{\upsilon}_{i}^{\mathcal{B}}\right)^{\top}\right\} \nonumber \\
& =\frac{1}{4}\sum_{i=1}^{n}s_{i}\left(1-\left(\hat{\upsilon}_{i}^{\mathcal{B}}\right)^{\top}\upsilon_{i}^{\mathcal{B}}\right)\label{eq:SO3PPF_RI_VM}
\end{align}
Let us introduce the following variable
\begin{align}
& \boldsymbol{\Upsilon}\left(M^{\mathcal{B}},\tilde{R}\right)\nonumber \\
& \hspace{2em}={\rm Tr}\left\{ \left(M^{\mathcal{B}}\right)^{-1}M^{\mathcal{B}}\tilde{R}\right\} \nonumber \\
& \hspace{2em}={\rm Tr}\left\{ \left(\sum_{i=1}^{n}s_{i}\upsilon_{i}^{\mathcal{B}}\left(\upsilon_{i}^{\mathcal{B}}\right)^{\top}\right)^{-1}\sum_{i=1}^{n}s_{i}\upsilon_{i}^{\mathcal{B}}\left(\hat{\upsilon}_{i}^{\mathcal{B}}\right)^{\top}\right\} \label{eq:SO3PPF_Gamma}
\end{align}
Consequently, any $\mathbf{vex}\left(\boldsymbol{\mathcal{P}}_{a}\left(M^{\mathcal{B}}\tilde{R}\right)\right)$,
$||M^{\mathcal{B}}\tilde{R}||_{I}$ and $\boldsymbol{\Upsilon}\left(M^{\mathcal{B}},\tilde{R}\right)$
will be obtained via a set of vectorial measurements as given in \eqref{eq:SO3PPF_VEX_VM},
\eqref{eq:SO3PPF_RI_VM}, and \eqref{eq:SO3PPF_Gamma}, respectively,
in all the subsequent calculations and derivations. Let us define
the minimum singular value of $\bar{\mathbf{M}}^{\mathcal{B}}$ as
$\underline{\lambda}:=\underline{\lambda}\left(\bar{\mathbf{M}}^{\mathcal{B}}\right)$,
$\mathcal{E}:=\mathcal{E}\left(||M^{\mathcal{B}}\tilde{R}||_{I},\xi\right)$,
and $\mu:=\mu\left(||M^{\mathcal{B}}\tilde{R}||_{I},\xi\right)$,
and consider the following filter kinematics
\begin{align}
\dot{\hat{R}}= & \hat{R}\left[\Omega_{m}-\hat{b}-W\right]_{\times},\quad\hat{R}\left(0\right)=\hat{R}_{0}\label{eq:SO3PPF_Rest_dot_VM}\\
\dot{\hat{b}}= & \frac{1}{2}\gamma\mu\mathcal{E}\mathbf{vex}\left(\boldsymbol{\mathcal{P}}_{a}\left(M^{\mathcal{B}}\tilde{R}\right)\right),\quad\hat{b}\left(0\right)=\hat{b}_{0}\label{eq:SO3PPF_b_est_VM}\\
W= & \frac{4}{\underline{\lambda}}\frac{k_{w}\mu\mathcal{E}-\dot{\xi}/\xi}{1+\boldsymbol{\Upsilon}\left(M^{\mathcal{B}},\tilde{R}\right)}\mathbf{vex}\left(\boldsymbol{\mathcal{P}}_{a}\left(M^{\mathcal{B}}\tilde{R}\right)\right)\label{eq:SO3PPF_W_VM}
\end{align}
where $\boldsymbol{\Upsilon}\left(M^{\mathcal{B}},\tilde{R}\right)$
and $\mathbf{vex}\left(\boldsymbol{\mathcal{P}}_{a}\left(M^{\mathcal{B}}\tilde{R}\right)\right)$
are defined in terms of vectorial measurements in \eqref{eq:SO3PPF_Gamma}
and \eqref{eq:SO3PPF_VEX_VM}, respectively, $\xi$ is a PPF defined
in \eqref{eq:SO3PPF_Presc}, $\mathcal{E}$ and $\mu$ are defined
in \eqref{eq:SO3PPF_trans3} and \eqref{eq:SO3PPF_Trans_Aux}, respectively,
with every $||\tilde{R}||_{I}$ being replaced by $||M^{\mathcal{B}}\tilde{R}||_{I}$,
$k_{w}$ and $\gamma$ are positive constants, and $\hat{b}$ is the
estimate of $b$. 
\begin{thm}
	\textbf{\label{thm:SO3PPF_2}} Consider the filter in \eqref{eq:SO3PPF_Rest_dot_VM},
	\eqref{eq:SO3PPF_b_est_VM} and \eqref{eq:SO3PPF_W_VM} to be coupled
	with the normalized vector measurements in \eqref{eq:SO3PPF_Vector_norm}
	and angular velocity measurements in \eqref{eq:SO3PPF_Angular} with
	the assumption that no noise is associated with the measurement $\Omega_{m}=\Omega+b$.
	Let two or more body-frame non-collinear vectors be available for
	measurements such that $M^{\mathcal{B}}$ is nonsingular. Define $\mathcal{U}\subseteq\mathbb{SO}\left(3\right)\times\mathbb{R}^{3}$
	by $\mathcal{U}:=\left\{ \left.\left(\tilde{R}_{0},\tilde{b}_{0}\right)\right|{\rm Tr}\left\{ \tilde{R}_{0}\right\} =-1,\tilde{b}_{0}=\underline{\boldsymbol{0}}_{3}\right\} $
	with $\tilde{R}_{0}=\tilde{R}\left(0\right)$ and $\tilde{b}_{0}=\tilde{b}\left(0\right)$.
	If $\tilde{R}_{0}\notin\mathcal{U}$ and $\mathcal{E}\left(0\right)\in\mathcal{L}_{\infty}$,
	then, all error signals are bounded, while $\mathcal{E}\left(t\right)$
	asymptotically approaches $0$ and $\tilde{R}$ asymptotically approaches
	$\mathbf{I}_{3}$.
\end{thm}
The observer dynamics in \eqref{eq:SO3PPF_Rest_dot_VM}, \eqref{eq:SO3PPF_b_est_VM}
and \eqref{eq:SO3PPF_W_VM} are guaranteed by Theorem \ref{thm:SO3PPF_2}
to be stable as $\mathcal{E}\left(t\right)$ approaches the origin
asymptotically. It follows that $\mathcal{E}\left(t\right)$ is bounded,
which in turn causes $||\tilde{R}||_{I}$ to obey the prescribed transient
and steady-state performance as described in \eqref{eq:SO3PPF_Presc}
in consistence with Remark \ref{rem:SO3PPF_1}.

\textbf{Proof. }Consider the error in attitude and bias being defined
similar to \eqref{eq:SO3PPF_R_error} and \eqref{eq:SO3PPF_b_tilde},
respectively. From \eqref{eq:SO3PPF_R_dynam} and \eqref{eq:SO3PPF_Rest_dot_Ry},
the error dynamics can be found to be analogous to \eqref{eq:SO3PPF_Rtilde_dot}.
From \eqref{eq:SO3PPF_MB_MI}, one can find the derivative of $M^{\mathcal{B}}$
to be
\begin{align}
\dot{M}^{\mathcal{B}} & =\dot{R}^{\top}M^{\mathcal{I}}R+R^{\top}M^{\mathcal{I}}\dot{R}\nonumber \\
& =-\left[\Omega\right]_{\times}R^{\top}M^{\mathcal{I}}R+R^{\top}M^{\mathcal{I}}R\left[\Omega\right]_{\times}\nonumber \\
& =-\left[\Omega\right]_{\times}M^{\mathcal{B}}+M^{\mathcal{B}}\left[\Omega\right]_{\times}\label{eq:SO3PPF_MB_dot}
\end{align}
Therefore, from \eqref{eq:SO3PPF_Rtilde_dot} and \eqref{eq:SO3PPF_MB_dot},
the derivative of $||M^{\mathcal{B}}\tilde{R}||_{I}$ can be expressed
as
\begin{align}
\frac{d}{dt}||M^{\mathcal{B}}\tilde{R}||_{I}= & -\frac{1}{4}{\rm Tr}\left\{ M^{\mathcal{B}}\dot{\tilde{R}}+\dot{M}^{\mathcal{B}}\tilde{R}\right\} \nonumber \\
= & -\frac{1}{4}{\rm Tr}\left\{ M^{\mathcal{B}}\left(\left[\tilde{R},\left[\Omega\right]_{\times}\right]+\tilde{R}\left[\tilde{b}-W\right]_{\times}\right)\right\} \nonumber \\
& -\frac{1}{4}{\rm Tr}\left\{ \left(-\left[\Omega\right]_{\times}M^{\mathcal{B}}+M^{\mathcal{B}}\left[\Omega\right]_{\times}\right)\tilde{R}\right\} \nonumber \\
= & -\frac{1}{4}{\rm Tr}\left\{ M^{\mathcal{B}}\tilde{R}\left[\tilde{b}-W\right]_{\times}\right\} \nonumber \\
& -\frac{1}{4}{\rm Tr}\left\{ \left[M^{\mathcal{B}}\tilde{R},\left[\Omega\right]_{\times}\right]\right\} \nonumber \\
= & \frac{1}{2}\mathbf{vex}\left(\boldsymbol{\mathcal{P}}_{a}\left(M^{\mathcal{B}}\tilde{R}\right)\right)^{\top}\left(\tilde{b}-W\right)\label{eq:SO3PPF_NormMBRtilde_dot}
\end{align}
where ${\rm Tr}\left\{ M^{\mathcal{B}}\tilde{R}\left[\tilde{b}\right]_{\times}\right\} =-2\mathbf{vex}\left(\boldsymbol{\mathcal{P}}_{a}\left(M^{\mathcal{B}}\tilde{R}\right)\right)^{\top}\tilde{b}$
as given in \eqref{eq:SO3PPF_Identity7}, and ${\rm Tr}\left\{ \left[M^{\mathcal{B}}\tilde{R},\left[\Omega\right]_{\times}\right]\right\} =0$
as defined in \eqref{eq:SO3PPF_Identity5}. Thus, in view of \eqref{eq:SO3PPF_NormR_dynam}
and \eqref{eq:SO3PPF_Trans_dot}, the derivative of the transformed
error in the sense of \eqref{eq:SO3PPF_NormRtilde_dot} can be found
to be
\begin{align}
\dot{\mathcal{E}}= & \frac{\mu}{2}\mathbf{vex}\left(\boldsymbol{\mathcal{P}}_{a}\left(M^{\mathcal{B}}\tilde{R}\right)\right)^{\top}\left(\tilde{b}-W\right)-\mu\frac{\dot{\xi}}{\xi}||M^{\mathcal{B}}\tilde{R}||_{I}\label{eq:SO3PPF_Trans_dot_VM}
\end{align}
Define the following candidate Lyapunov function as
\begin{align}
V\left(\mathcal{E},\tilde{b}\right) & =\frac{1}{2}\mathcal{E}^{2}+\frac{1}{2\gamma}||\tilde{b}||^{2}\label{eq:SO3PPF_V_VM}
\end{align}
The derivative of $V:=V\left(\mathcal{E},\tilde{b}\right)$ in \eqref{eq:SO3PPF_V_VM}
can be expressed as
\begin{align}
\dot{V}= & \mathcal{E}\dot{\mathcal{E}}-\frac{1}{\gamma}\tilde{b}^{\top}\dot{\hat{b}}\nonumber \\
= & \mathcal{E}\mu\left(\frac{1}{2}\mathbf{vex}\left(\boldsymbol{\mathcal{P}}_{a}\left(M^{\mathcal{B}}\tilde{R}\right)\right)^{\top}\left(\tilde{b}-W\right)-\frac{\dot{\xi}}{\xi}||M^{\mathcal{B}}\tilde{R}||_{I}\right)\nonumber \\
& -\frac{1}{\gamma}\tilde{b}^{\top}\dot{\hat{b}}\label{eq:SO3PPF_Vdot_VM-1}
\end{align}
Directly substituting for $\dot{\hat{b}}$ and $W$ in \eqref{eq:SO3PPF_b_est_VM},
and \eqref{eq:SO3PPF_W_VM}, respectively, one obtains
\begin{align}
\dot{V}\leq & \frac{\dot{\xi}}{\xi}\left(\frac{2}{\underline{\lambda}}\frac{\left\Vert \mathbf{vex}\left(\boldsymbol{\mathcal{P}}_{a}\left(M^{\mathcal{B}}\tilde{R}\right)\right)\right\Vert ^{2}}{1+\boldsymbol{\Upsilon}\left(M^{\mathcal{B}},\tilde{R}\right)}-\left\Vert M^{\mathcal{B}}\tilde{R}\right\Vert _{I}\right)\mu\mathcal{E}\nonumber \\
& -\frac{2}{\underline{\lambda}}\frac{k_{w}\mu^{2}\mathcal{E}^{2}}{1+\boldsymbol{\Upsilon}\left(M^{\mathcal{B}},\tilde{R}\right)}\left\Vert \mathbf{vex}\left(\boldsymbol{\mathcal{P}}_{a}\left(M^{\mathcal{B}}\tilde{R}\right)\right)\right\Vert ^{2}\label{eq:SO3PPF_Vdot_VM-2}
\end{align}
One can also easily find
\begin{equation}
\frac{\dot{\xi}}{\xi}\left(\frac{2}{\underline{\lambda}}\frac{\left\Vert \mathbf{vex}\left(\boldsymbol{\mathcal{P}}_{a}\left(M^{\mathcal{B}}\tilde{R}\right)\right)\right\Vert ^{2}}{1+\boldsymbol{\Upsilon}\left(M^{\mathcal{B}},\tilde{R}\right)}-\left\Vert M^{\mathcal{B}}\tilde{R}\right\Vert _{I}\right)\mu\mathcal{E}\leq0\label{eq:SO3PPF_Factor}
\end{equation}
where $\mathcal{E}>0\forall||M^{\mathcal{B}}\tilde{R}||_{I}\neq0$
and $\mathcal{E}=0$ at $||M^{\mathcal{B}}\tilde{R}||_{I}=0$ as given
in (i) Proposition \ref{Prop:SO3PPF_1}, and $\mu>0\forall t\geq0$
as given in \eqref{eq:SO3PPF_Trans_Aux}. Also, $\dot{\xi}$ is a
negative strictly increasing component which satisfies $\dot{\xi}\rightarrow0$
as $t\rightarrow\infty$, and $\xi:\mathbb{R}_{+}\to\mathbb{R}_{+}$
such that $\xi\rightarrow\xi_{\infty}$ as $t\rightarrow\infty$.
Thus, $\dot{\xi}/\xi\leq0$. In addition, consider \eqref{eq:SO3PPF_lemm1_3}
in Lemma \ref{Lemm:SO3PPF_1}, the expression in \eqref{eq:SO3PPF_Factor}
is negative semi-definite. Consequently, the inequality in \eqref{eq:SO3PPF_Vdot_VM-2}
can be expressed as
\begin{align}
\dot{V}\leq & -k_{w}\mu^{2}\mathcal{E}^{2}\left\Vert M^{\mathcal{B}}\tilde{R}\right\Vert _{I}\label{eq:SO3PPF_Vdot_VM_Final}
\end{align}
This implies that $V\left(t\right)\leq V\left(0\right),\forall t\geq0$.
Given that $\tilde{R}_{0}\notin\mathcal{U}$, $\tilde{b}$ is bounded
for $t\geq0$, and $\mathcal{E}\in\mathcal{L}_{\infty},\forall t\geq0$.
As such, $\mathcal{E}$ remains bounded and well-defined for all $t\geq0$.
In order to prove asymptotic convergence of $\mathcal{E}$ to the
origin and $\tilde{R}$ to the identity for all $\tilde{R}_{0}\notin\mathcal{U}$,
one obtains the second derivative of \eqref{eq:SO3PPF_V_VM} as 
\begin{align}
\ddot{V}\leq & -2k_{w}\left(\mathcal{E}\dot{\mathcal{E}}\mu^{2}+\mathcal{E}^{2}\mu\dot{\mu}\right)||M^{\mathcal{B}}\tilde{R}||_{I}\nonumber \\
& -k_{w}\mathcal{E}^{2}\mu^{2}\frac{d}{dt}||M^{\mathcal{B}}\tilde{R}||_{I}\label{eq:SO3PPF_Vddot_VM}
\end{align}
Consider the result in \eqref{eq:SO3PPF_Trans_Aux}, as such, it can
be shown that
\begin{equation}
\dot{\mu}=-\frac{1}{2}\frac{\underline{\delta}\dot{\xi}+\frac{d}{dt}||M^{\mathcal{B}}\tilde{R}||_{I}}{\left(\underline{\delta}\xi+||\tilde{R}||_{I}\right)^{2}}-\frac{1}{2}\frac{\bar{\delta}\dot{\xi}-\frac{d}{dt}||M^{\mathcal{B}}\tilde{R}||_{I}}{\left(\bar{\delta}\xi-||\tilde{R}||_{I}\right)^{2}}\label{eq:SO3PPF_mu_VMdot}
\end{equation}
with $\dot{\xi}=-\ell\left(\xi^{0}-\xi^{\infty}\right)\exp\left(-\ell t\right)$.
Due to the fact that $||\dot{\tilde{R}}||_{I}$ is bounded, $\dot{\mu}$
is bounded and in turn $\ddot{V}$ is bounded for all $t\geq0$. Thus,
$\dot{V}$ is uniformly continuous and in accordance with Barbalat
Lemma, $\dot{V}\rightarrow0$ implies that either $||\mathcal{E}||\rightarrow0$
or $||M^{\mathcal{B}}\tilde{R}||_{I}\rightarrow0$ or both $||\mathcal{E}||\rightarrow0$
and $||M^{\mathcal{B}}\tilde{R}||_{I}\rightarrow0$ as $t\rightarrow\infty$.
From property (i) and (ii) of Proposition \ref{Prop:SO3PPF_1}, $||\mathcal{E}||\rightarrow0$
indicates that $||M^{\mathcal{B}}\tilde{R}||_{I}\rightarrow0$ and
vice versa. Thus, $\dot{V}\rightarrow0$ implies that $||\mathcal{E}||\rightarrow0$
and $||M^{\mathcal{B}}\tilde{R}||_{I}\rightarrow0$, which means that
$\tilde{R}$ asymptotically approaches $\mathbf{I}_{3}$ consistent
with property (iii) of Proposition \ref{Prop:SO3PPF_1}, which completes
the proof.

It is clear that the gains associated with the vex operator of $\dot{\hat{b}}$
and $W$ in \eqref{eq:SO3PPF_b_est_Ry}, and \eqref{eq:SO3PPF_W_Ry},
or in \eqref{eq:SO3PPF_b_est_VM}, and \eqref{eq:SO3PPF_W_VM}, respectively,
are dynamic. Their values rely on $\mu$, $\mathcal{E}$ and $||\tilde{R}||_{I}$
or $||M^{\mathcal{B}}\tilde{R}||_{I}$. Their dynamic behavior has
the essential role of forcing the proposed observer to comply with
the prescribed performance constraints. Thus, the proposed filter
has a remarkable advantage which is reflected in the dynamic gains
becoming increasingly aggressive as $||\tilde{R}||_{I}$ approaches
the unstable equilibria $+1$. On the other side, these gains reduce
significantly as $\mathcal{E}\rightarrow0$. These dynamic gains directly
impact the proposed nonlinear filter forcing it to adhere to the predefined
prescribed performance features imposed by the user and thereby satisfying
the predefined measures of transient as well as steady-state measures.
\begin{rem}
	\label{rem:SO3PPF_3}\textbf{(Notes on filter design parameters)}
	$\bar{\delta}$, $\underline{\delta}$, and $\xi_{0}$ define the
	dynamic boundaries of the transformed error $\mathcal{E}$. $\xi_{0}$
	and $\xi_{\infty}$ refer to the boundaries of the large and small
	sets, respectively. $\ell$ controls the convergence rate of the dynamic
	boundaries from large to narrow set. The asymptotic convergence of
	$||\tilde{R}||_{I}$ or $||M^{\mathcal{B}}\tilde{R}||_{I}$ is guaranteed
	by selecting $\bar{\delta}=\underline{\delta}$. Also, increasing
	the value of $\ell$ would lead to faster rate of convergence of $||\tilde{R}||_{I}$
	or $||M^{\mathcal{B}}\tilde{R}||_{I}$ to the origin. It should be
	noted that if the initial value of $||\tilde{R}\left(0\right)||_{I}$
	or $||M^{\mathcal{B}}\tilde{R}\left(0\right)||_{I}$ are unknown,
	the user could select $\bar{\delta}$, $\underline{\delta}$, and
	$\xi_{0}$ based on the highest value of $||\tilde{R}\left(0\right)||_{I}$,
	therefore accounting for the worst possible scenario, since $||\tilde{R}\left(0\right)||_{I}\in\left[0,1\right]$,
	and thus the prescribed performance is guaranteed.
\end{rem}
The filter design algorithm proposed in Subsection \ref{subsec:SO3PPF_Explicit-Filter}
can be summarized briefly as 
\begin{enumerate}
	\item[\textbf{A.1}]  Select $\bar{\delta}=\underline{\delta}>||M^{\mathcal{B}}\tilde{R}\left(0\right)||_{I}$,
	the ultimate bound of the small set of the desired steady-state error
	$\xi_{\infty}$ and the desired convergence rate $\ell$. 
	\item[\textbf{A.2}]  Evaluate the vex operator $\mathbf{vex}\left(\boldsymbol{\mathcal{P}}_{a}\left(M^{\mathcal{B}}\tilde{R}\right)\right)$,
	the normalized Euclidean distance error $||M^{\mathcal{B}}\tilde{R}||_{I}$,
	and $\boldsymbol{\Upsilon}\left(M^{\mathcal{B}},\tilde{R}\right)$
	from \eqref{eq:SO3PPF_VEX_VM}, \eqref{eq:SO3PPF_RI_VM}, and \eqref{eq:SO3PPF_Gamma},
	respectively, in the form of vector measurements. 
	\item[\textbf{A.3}]  Evaluate the prescribed performance function $\xi$ from equation
	\eqref{eq:SO3PPF_Presc}. 
	\item[\textbf{A.4}]  Evaluate $\mu\left(||M^{\mathcal{B}}\tilde{R}||_{I},\xi\right)$
	and $\mathcal{E}\left(||M^{\mathcal{B}}\tilde{R}||_{I},\xi\right)$
	from equations \eqref{eq:SO3PPF_Trans_Aux} and \eqref{eq:SO3PPF_trans3},
	respectively. 
	\item[\textbf{A.5}]  Evaluate the filter design $\dot{\hat{R}}$, $\dot{\hat{b}}$ and
	$W$ from \eqref{eq:SO3PPF_Rest_dot_VM}, \eqref{eq:SO3PPF_b_est_VM},
	and \eqref{eq:SO3PPF_W_VM}, respectively. 
	\item[\textbf{A.6}]  Go to \textbf{A.2}. 
\end{enumerate}
The same steps can be applied for the filter design in Subsection
\ref{subsec:SO3PPF_Passive-Filter}.

\section{Simulations \label{sec:SO3PPF_Simulations}}

The performance of the two proposed nonlinear attitude filters on
$\mathbb{SO}\left(3\right)$ with predefined measures is presented
in this section considering large error initialization and high level
of noise and bias in the measurements. In this regard, consider the
set of measurements given as follows:
\[
\begin{cases}
{\rm v}_{i}^{\mathcal{B}} & =R^{\top}{\rm v}_{i}^{\mathcal{I}}+{\rm b}_{i}^{\mathcal{B}}+\omega_{i}^{\mathcal{B}}\\
\Omega_{m} & =\Omega+b+\omega
\end{cases}
\]
which exemplifies a set measurements obtained from a low-cost IMUs
module, for all $i=1,2$. Let the rotational matrix $R$ be acquired
from attitude dynamics in equation \eqref{eq:SO3PPF_R_dynam} and
suppose that the input signal of the angular velocity is given by
\[
\Omega=\left[\begin{array}{c}
{\rm sin}\left(0.7t\right)\\
0.7{\rm sin}\left(0.5t+\pi\right)\\
0.5{\rm sin}\left(0.3t+\frac{\pi}{3}\right)
\end{array}\right]\left({\rm rad/sec}\right)
\]
with $R\left(0\right)=\mathbf{I}_{3}$ being the initial attitude.
Consider that a wide-band of a zero mean random noise process vector
with standard deviation (STD) of $0.2\left({\rm rad/sec}\right)$
and bias $b=0.1\left[1,-1,1\right]^{\top}$ is contaminating the true
angular velocity $\left(\Omega\right)$ such that $\Omega_{m}=\Omega+b+\omega$.
Let two non-collinear inertial frame vectors be given by ${\rm v}_{1}^{\mathcal{I}}=\frac{1}{\sqrt{3}}\left[1,-1,1\right]^{\top}$
and ${\rm v}_{2}^{\mathcal{I}}=\left[0,0,1\right]^{\top}$, whereas
the body-frame vectors ${\rm v}_{1}^{\mathcal{B}}$ and ${\rm v}_{2}^{\mathcal{B}}$
are given by ${\rm v}_{i}^{\mathcal{B}}=R^{\top}{\rm v}_{i}^{\mathcal{I}}+{\rm b}_{i}^{\mathcal{B}}+\omega_{i}^{\mathcal{B}}$
for all $i=1,2$. Similarly, suppose that an additional zero mean
random noise vector $\omega_{i}^{\mathcal{B}}$ with ${\rm STD=}0.08$
corrupts the body-frame vector measurements with bias components ${\rm b}_{1}^{\mathcal{B}}=0.1\left[-1,1,0.5\right]^{\top}$
and ${\rm b}_{2}^{\mathcal{B}}=0.1\left[0,0,1\right]^{\top}$. ${\rm v}_{i}^{\mathcal{I}}$
and ${\rm v}_{i}^{\mathcal{B}}$ are normalized and the third vector
is extracted by $\upsilon_{3}^{\mathcal{I}}=\upsilon_{1}^{\mathcal{I}}\times\upsilon_{2}^{\mathcal{I}}$
and $\upsilon_{3}^{\mathcal{B}}=\upsilon_{1}^{\mathcal{B}}\times\upsilon_{2}^{\mathcal{B}}$.
The confidence level of body-frame measurements was chosen as $s_{1}=1.4$,
$s_{2}=1.4$, and $s_{3}=0.2$. For the semi-direct filter in Subsection
\ref{subsec:SO3PPF_Passive-Filter}, the corrupted reconstructed attitude
$R_{y}$ is defined using SVD \nameref{sec:SO3PPF_AppendixB} or see
the Appendix in \cite{hashim2018SO3Stochastic, hashim2018SE3Stochastic} where $\tilde{R}=R_{y}^{\top}\hat{R}$.

To illustrate the robustness of the proposed filtering algorithms,
a very large initial attitude error is considered. The initial rotation
of the attitude estimate is defined in accordance with angle-axis
parameterization in \eqref{eq:SO3PPF_att_ang} as $\hat{R}\left(0\right)=\mathcal{R}_{\alpha}\left(\alpha,u/||u||\right)$
with $\alpha=178\left({\rm deg}\right)$ and $u$= $\left[4,1,5\right]^{\top}$.
As such, $||\tilde{R}||_{I}\approx0.9999$ which is very close to
the unstable equilibria. Initial bias estimate is $\hat{b}\left(0\right)=\left[0,0,0\right]^{\top}$.
The design parameters are chosen as $\gamma=1$, $k_{w}=3$, $\bar{\delta}=\underline{\delta}=1.2$,
$\xi_{0}=1.2$, $\xi_{\infty}=0.05$, and $\ell=3$. The total time
of the simulation is 15 seconds.

The color notation is as follows: green color represents a true value,
red depicts the performance of the nonlinear semi-direct filter on
$\mathbb{SO}\left(3\right)$ derived using a group of vectorial measurements
and reconstructed attitude as described in Subsection \ref{subsec:SO3PPF_Passive-Filter},
and blue demonstrates the performance of the direct filter characterized
in Subsection \ref{subsec:SO3PPF_Explicit-Filter} which does not
demand attitude reconstruction. Also, magenta describes a measured
value while orange and purple refer to prescribed performance response.

\begin{figure}[h!]
	\centering{}\includegraphics[scale=0.26]{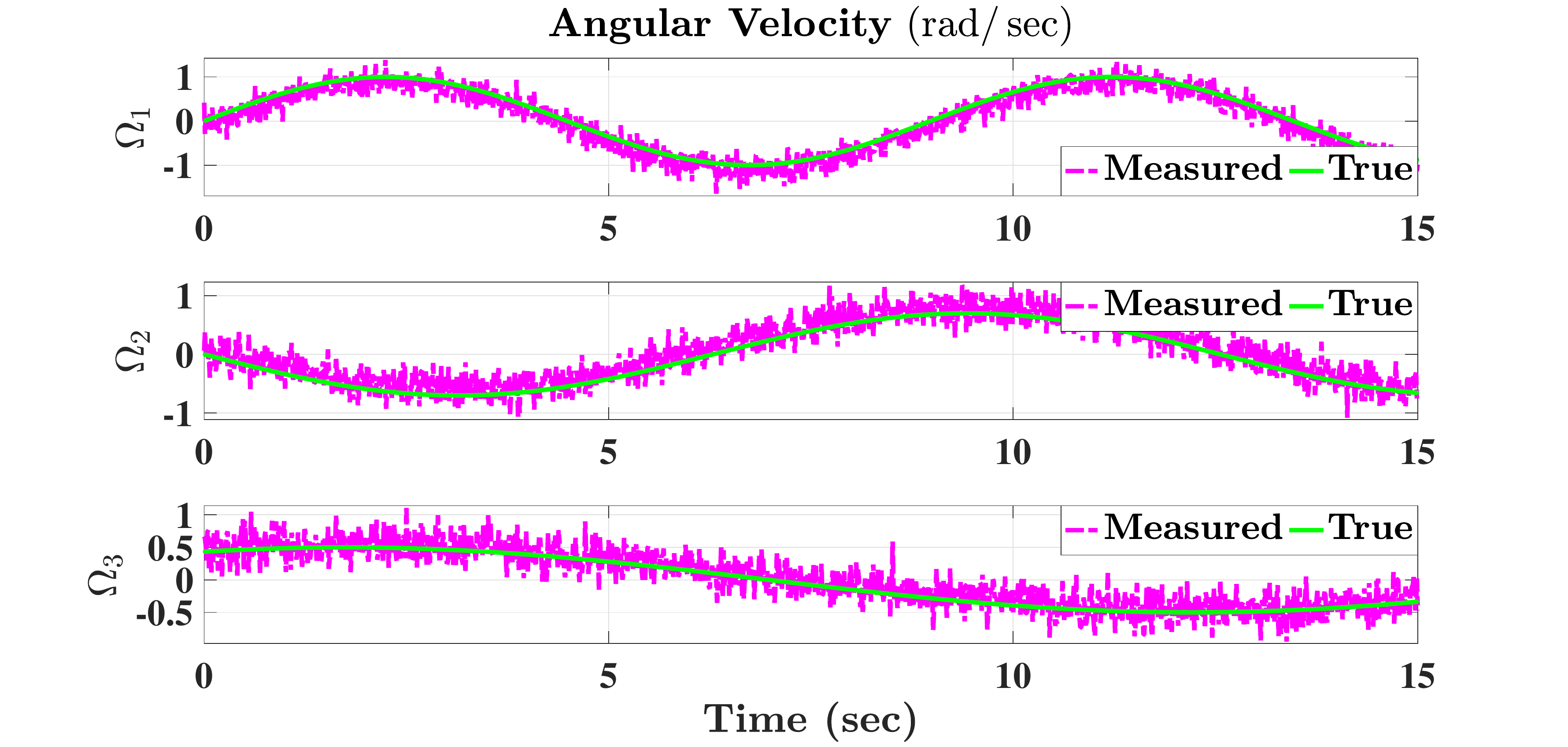}\caption{True and measured angular velocities.}
	\label{fig:SO3PPF_3} 
\end{figure}

\begin{figure}[h!]
	\centering{}\includegraphics[scale=0.26]{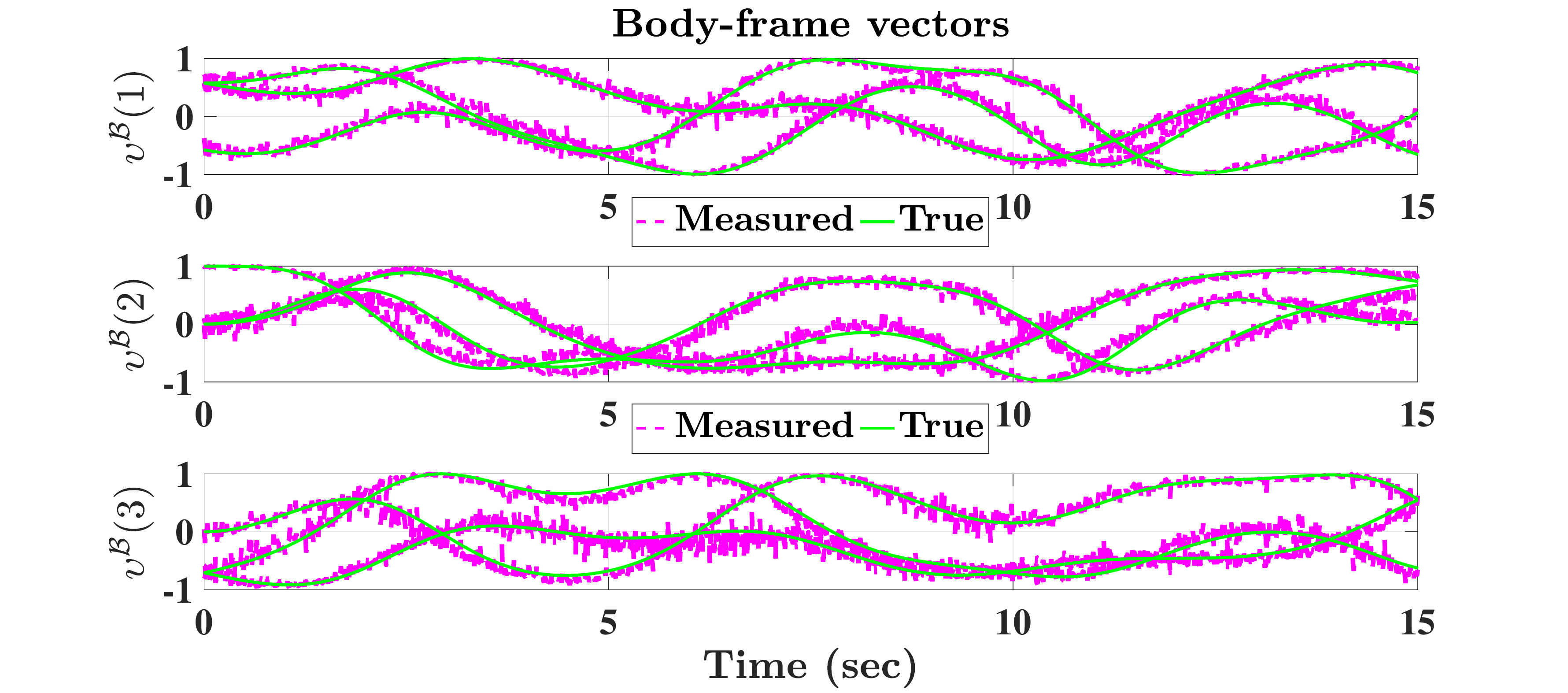}\caption{Body-frame vectorial measurements: true and measured. }
	\label{fig:SO3PPF_4} 
\end{figure}

\begin{figure*}[h!]
	\centering{}\includegraphics[scale=0.42]{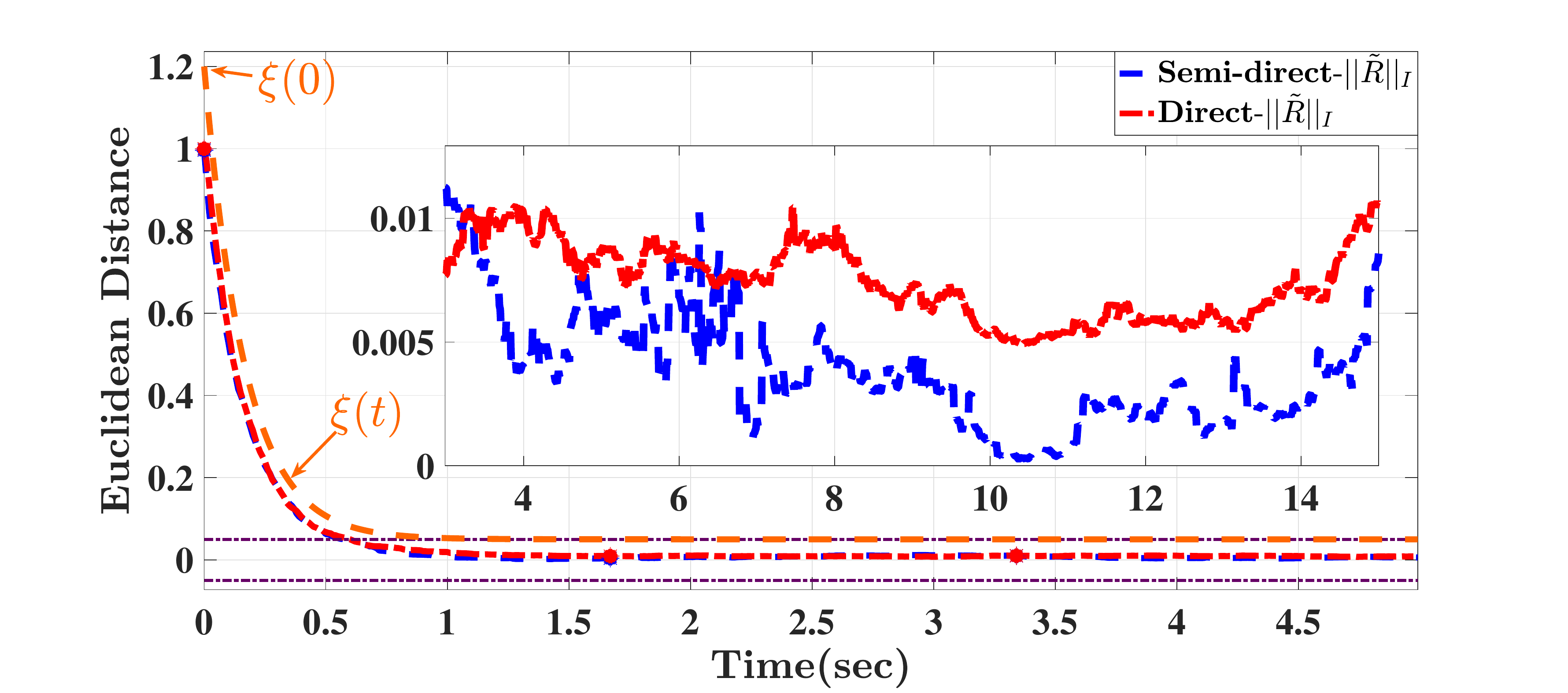}\caption{Transient and steady-state performance of normalized Euclidean distance.}
	\label{fig:SO3PPF_5} 
\end{figure*}

\begin{figure*}[h!]
	\centering{}\includegraphics[scale=0.38]{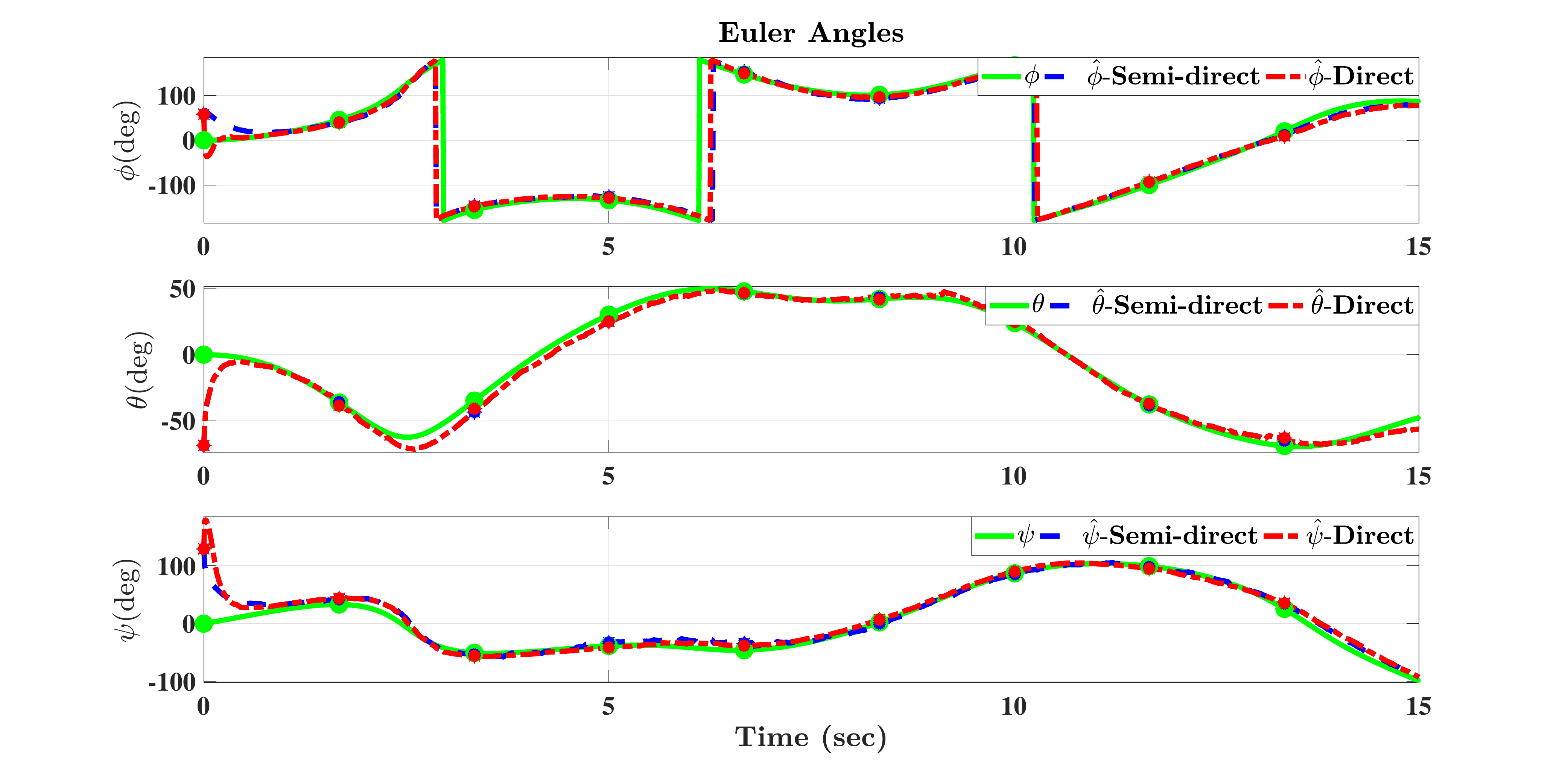}\caption{Tracking performance of Euler angles (roll $(\phi)$ and pitch $(\theta)$,
		yaw $(\psi)$).}
	\label{fig:SO3PPF_6} 
\end{figure*}

\begin{figure}[h!]
	\centering{}\includegraphics[scale=0.26]{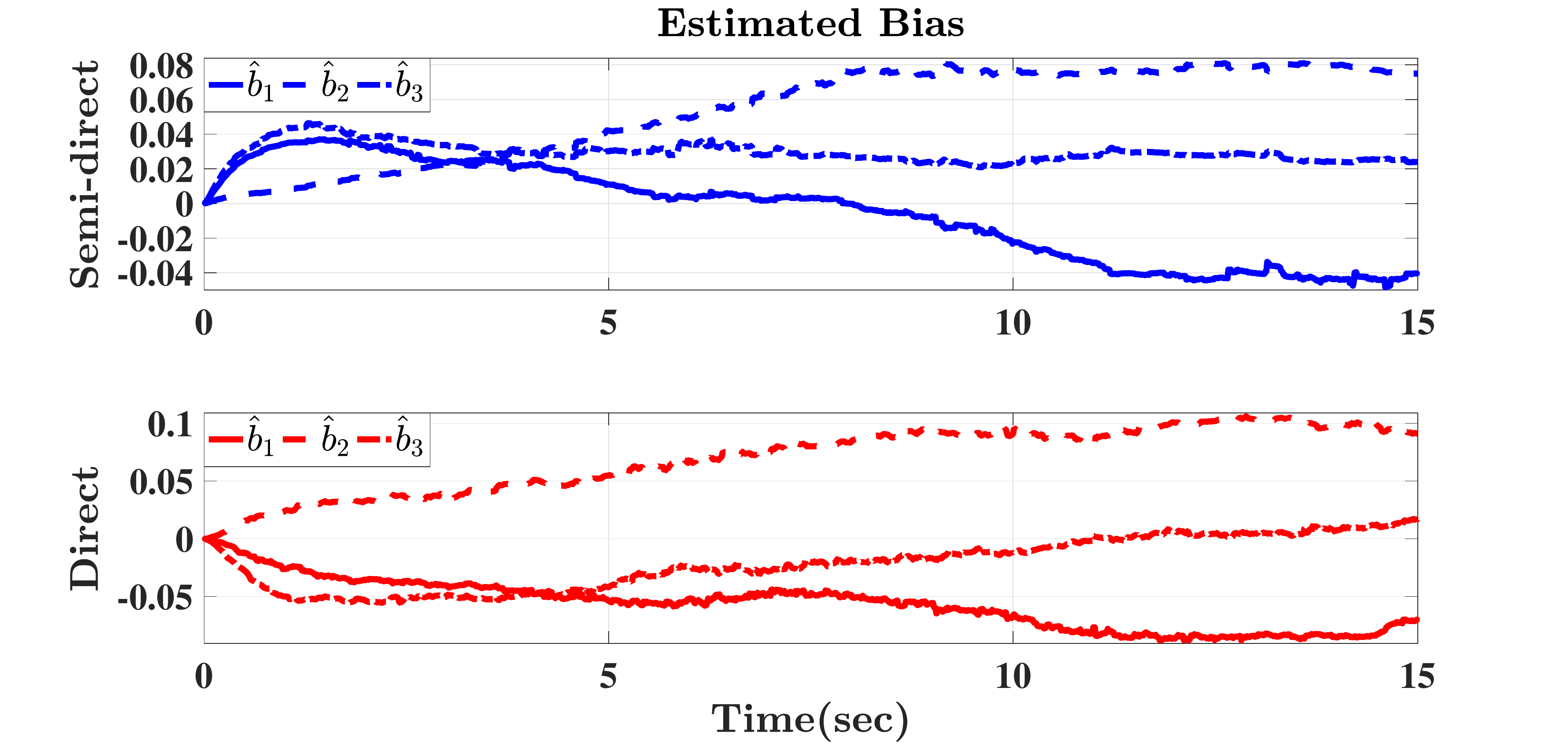}\caption{The estimated bias of the proposed filters.}
	\label{fig:SO3PPF_7} 
\end{figure}

Figure \ref{fig:SO3PPF_3} and \ref{fig:SO3PPF_4} illustrate high
values of noise and bias components present in angular velocity and
body-frame vector measurements plotted against the true values. Figure
\ref{fig:SO3PPF_5} illustrates the systematic and smooth convergence
of the normalized Euclidean distance error $||\tilde{R}||_{I}$. It
can be noticed in Figure \ref{fig:SO3PPF_5} that the error function
for $||\tilde{R}||_{I}=\frac{1}{4}{\rm Tr}\left\{ \mathbf{I}_{3}-R^{\top}\hat{R}\right\} $
started very near to the unstable equilibria within a given large
set and ended within a given small residual set obeying the PPF. Thus,
Figure \ref{fig:SO3PPF_5} confirms the stability analysis discussed
in the previous section and illustrates the robustness of the proposed
filter. The output performance of the proposed filters in Euler angles
representation is shown in Figure \ref{fig:SO3PPF_6}. The three Euler
angles $\left(\phi,\theta,\psi\right)$ in Figure \ref{fig:SO3PPF_6}
show impressive tracking performance with fast convergence to the
true angles. Finally, the boundedness of the estimated bias $\hat{b}$
is illustrated in Figure \ref{fig:SO3PPF_7}.

Table \ref{tab:SO3PPF_1} contains a synopsis of statistical details
of the mean and the STD of the error ($||\tilde{R}||_{I}$). These
details facilitate the comparison of the steady-state performance
of the two filters proposed in this paper with respect to $||\tilde{R}||_{I}$.
In spite of the fact that both filters have extremely small mean of
$||\tilde{R}||_{I}$, the semi-direct attitude filter with prescribed
performance showed a remarkably smaller mean errors and STD when compared
to the direct attitude filter with prescribed performance. Numerical
results outlined in Table \ref{tab:SO3PPF_1} demonstrate effectiveness
and robustness of the proposed nonlinear attitude filters against
large error initialization and uncertainties in sensor measurements
as illustrated in Figure \ref{fig:SO3PPF_3}, \ref{fig:SO3PPF_4},
\ref{fig:SO3PPF_6}, \ref{fig:SO3PPF_5}, and \ref{fig:SO3PPF_7}. 

\begin{table}[H]
	\caption{\label{tab:SO3PPF_1}Statistical analysis of $||\tilde{R}||_{I}$
		of the proposed two filters.}
	
	\centering{}%
	\begin{tabular}{c|>{\centering}p{2.5cm}|>{\centering}p{2.5cm}}
		\hline 
		\multicolumn{3}{c}{Output data of $||\tilde{R}||_{I}$ over the period (1-15 sec)}\tabularnewline
		\hline 
		\hline 
		Filter  & Semi-direct  & Direct\tabularnewline
		\hline 
		Mean  & $4.2\times10^{-3}$  & $6.9\times10^{-3}$\tabularnewline
		\hline 
		STD  & $2.5\times10^{-3}$  & $2.1\times10^{-3}$\tabularnewline
		\hline 
	\end{tabular}
\end{table}

The robustness and the superior convergence properties of the proposed
nonlinear attitude filters with guaranteed performance are presented
and compared to a well-known nonlinear attitude complimentary filter
termed nonlinear passive complementary filter \cite{mahony2008nonlinear}
as well as to a standard attitude filter which belongs to the family
of Gaussian attitude filters and is termed multiplicative extended
Kalman filter (MEKF) \cite{markley2003attitude} in Subsection \ref{subsec:Comp1}
and \ref{subsec:Comp2}, respectively.

\subsection{Proposed Filters vs Nonlinear Attitude Filters \label{subsec:Comp1}}

To further illustrate the robustness and the superior convergence
properties of the proposed nonlinear attitude filters as opposed to
the conventional nonlinear attitude filters, a fair comparison is
presented. Consider the following nonlinear passive complementary
filter given in \cite{mahony2008nonlinear}
\begin{equation}
\begin{cases}
\dot{\hat{R}} & =\hat{R}\left[\Omega_{m}-\hat{b}-W\right]_{\times},\quad\hat{R}\left(0\right)=\hat{R}_{0}\\
\dot{\hat{b}} & =k_{1}\mathbf{vex}\left(\boldsymbol{\mathcal{P}}_{a}\left(\tilde{R}\right)\right),\quad\hat{b}\left(0\right)=\hat{b}_{0},\tilde{R}=R_{y}^{\top}\hat{R}\\
W & =k_{1}\mathbf{vex}\left(\boldsymbol{\mathcal{P}}_{a}\left(\tilde{R}\right)\right),\tilde{R}=R_{y}^{\top}\hat{R}
\end{cases}\label{eq:Passive}
\end{equation}
where $k_{1}>0$. A fair comparison between the proposed semi-direct
attitude filter and the nonlinear passive complementary filter in
\cite{mahony2008nonlinear} is attainable due to the shared structure
of the filters. Consider initializing the nonlinear passive complementary
filter analogously to the semi-direct attitude filter given at the
beginning of the Simulation Section. To ensure validity of the comparison,
three variations of the design parameter $k_{1}$ in \eqref{eq:Passive}
namely, $k_{1}=1$, $k_{1}=10$ and $k_{1}=100$. 
\begin{figure*}[h]
	\centering{}\includegraphics[scale=0.45]{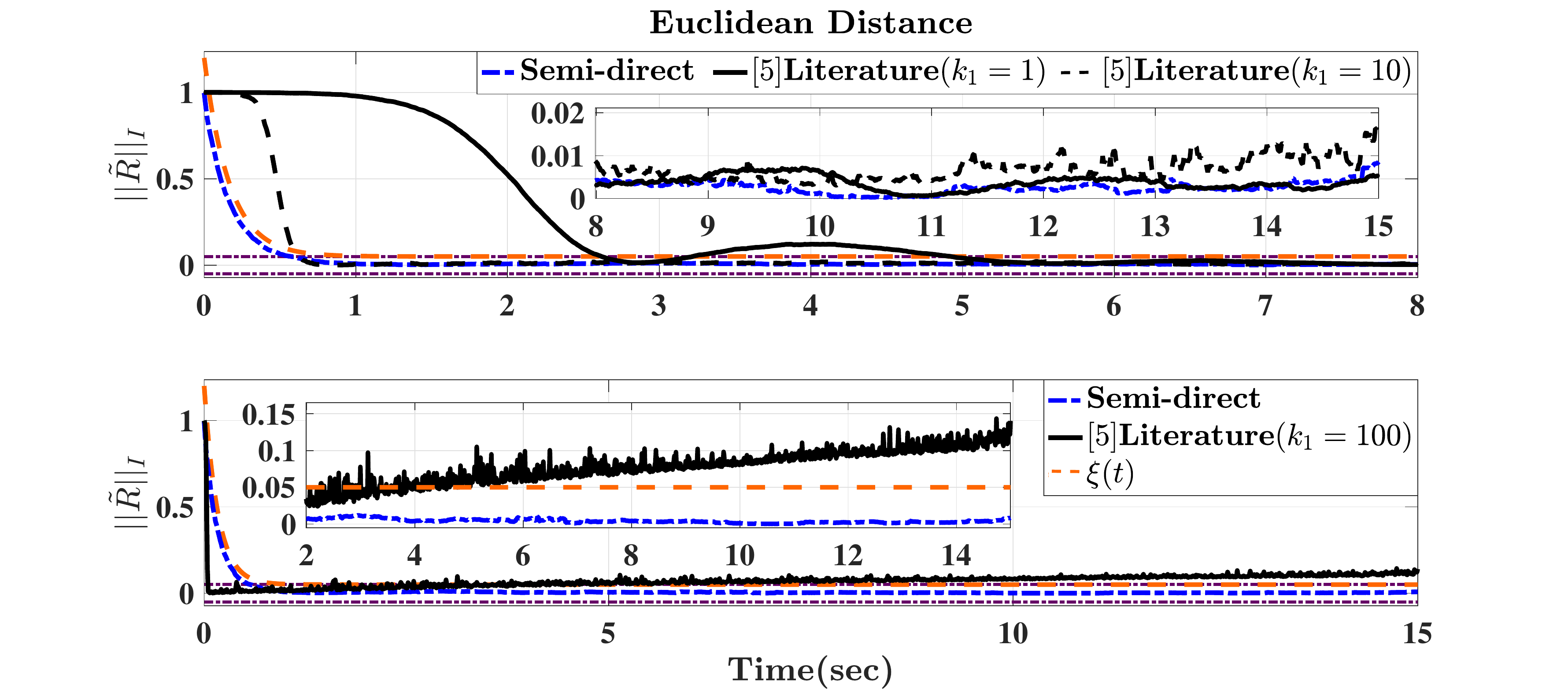}\caption{Transient and steady-state performance of normalized Euclidean distance:
		Semi-direct filter vs literature \cite{mahony2008nonlinear}.}
	\label{fig:SO3PPF_8} 
\end{figure*}
In this Subsection,
the color notation is as follows: black solid and dashed lines describe
the performance of the nonlinear passive complementary filter, blue
center line depicts the proposed semi-direct attitude filter while
orange and purple refer to the prescribed performance response. It
can be noticed in the upper portion of Figure \ref{fig:SO3PPF_8}
that smaller value of $k_{1}$ results in slower transient performance
with less oscillatory behavior in the steady-state. In contrast, the
lower portion of Figure \ref{fig:SO3PPF_8} illustrates that higher
value of $k_{1}$ leads to faster transient performance with higher
levels of oscillation in the steady-state. Moreover, Figure \ref{fig:SO3PPF_8}
shows that the predefined measure of transient performance cannot
be achieved for low value of $k_{1}$, since the transient performance
of the passive complementary filter violates the dynamic reducing
boundaries. In the same spirit, the predefined characteristics of
steady-state performance cannot be achieved for high value of $k_{1}$.
These results confirm Remark \ref{rem:RemNew}. 

Therefore, the nonlinear
attitude filters given in the literature, for example \cite{mahony2005complementary,mahony2008nonlinear,hamel2006attitude,grip2012attitude,lee2012exponential,zlotnik2017nonlinear}
cannot guarantee a predefined measure of convergence properties. The
semi-direct attitude filter, on the other side, obeys the dynamically
reducing boundaries and allows to achieve a desired level of prescribed
performance.

Table \ref{tab:SO3PPF_2} compares the statistical details, namely
the mean and the STD of $||\tilde{R}||_{I}$, of the proposed semi-direct
attitude filter and the nonlinear passive complementary filter. The
above-mentioned statistics describe the output performance with respect
to $||\tilde{R}||_{I}$ over the steady-state period of time depicted
in Figure \ref{fig:SO3PPF_8}. The semi-direct attitude filter displays
smaller values of mean and STD of $||\tilde{R}||_{I}$ when compared
to the passive complementary filter for all the considered cases of
$k_{1}=1$, $k_{1}=10$ and $k_{1}=100$. Moreover, the numerical
results listed in Table \ref{tab:SO3PPF_2} illustrate the effectiveness
and robustness of the proposed nonlinear attitude filters against
large error initialization and uncertainties in sensor measurements
which make them a good fit for measurements obtained from low-cost
IMUs modules.
\begin{table}[H]
	\caption{\label{tab:SO3PPF_2} Statistical analysis of $||\tilde{R}||_{I}$
		of the semi-direct filter vs literature.}
	
	\centering{}%
	\begin{tabular}{c|c|c|c|c}
		\hline 
		\multicolumn{5}{c}{Output data of $||\tilde{R}||_{I}$ over the period (7-15 sec)}\tabularnewline
		\hline 
		\hline 
		\multirow{2}{*}{Filter} & \multirow{2}{*}{Semi-direct} & \multicolumn{3}{c}{Passive Filter \cite{mahony2008nonlinear}}\tabularnewline
		\cline{3-5} 
		&  & $k_{1}=1$ & $k_{1}=10$ & $k_{1}=100$\tabularnewline
		\hline 
		Mean & $2.7\times10^{-3}$ & $4.5\times10^{-3}$ & $6.9\times10^{-3}$ & $91.9\times10^{-3}$\tabularnewline
		\hline 
		STD & $1.4\times10^{-3}$ & $2.9\times10^{-3}$ & $2.7\times10^{-3}$ & $14.2\times10^{-3}$\tabularnewline
		\hline 
	\end{tabular}
\end{table}

\subsection{Proposed Filters vs Gaussian Attitude Filters \label{subsec:Comp2}}

\begin{figure*}[h]
	\centering{}\includegraphics[scale=0.45]{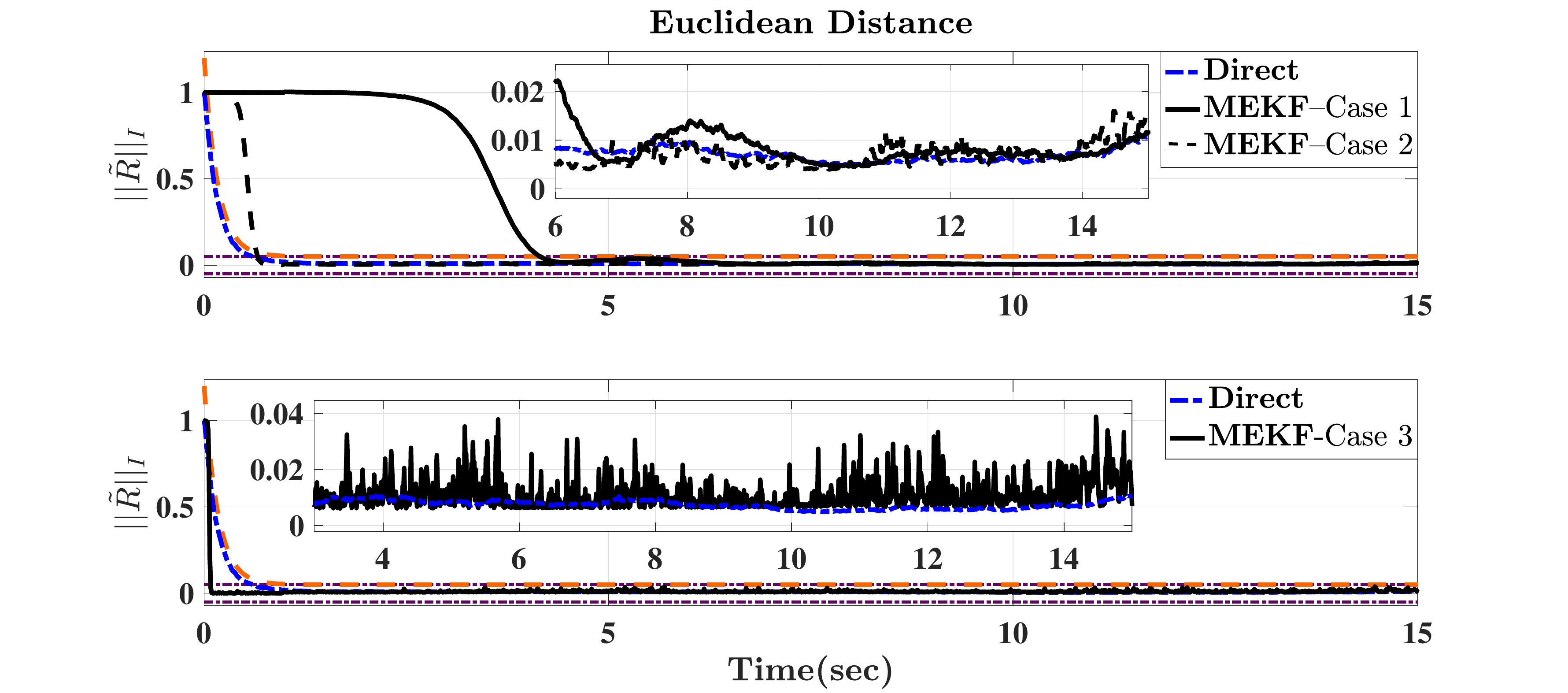}\caption{Transient and steady-state performance of normalized Euclidean distance:
		Direct filter vs MEKF \cite{markley2003attitude}.}
	\label{fig:SO3PPF_9} 
\end{figure*}

In this subsection the effectiveness and the high convergence capabilities
of the proposed nonlinear attitude filters are compared to the performance
of a Gaussian attitude filter. A comparison between the proposed direct
attitude filter and the MEKF in \nameref{sec:SO3PPF_AppendixC} is
presented. Consider the MEKF in \nameref{sec:SO3PPF_AppendixC} initialized
similar to the direct attitude filter given at the beginning of the
Simulation Section. To guarantee validity of the comparison, three
cases of the design parameters of MEKF have been detailed in Table
\ref{tab:SO3PPF_3}.
\begin{table}[H]
	\caption{\label{tab:SO3PPF_3} MEKF design parameters.}
	
	\noindent \centering{}%
	\begin{tabular}{c|c|c|c}
		\hline 
		Case & \multicolumn{3}{c}{Design Parameters}\tabularnewline
		\hline 
		\hline 
		Case 1 & $\mathcal{\bar{Q}}_{v\left(i\right)}=\mathbf{I}_{3}$ & $\mathcal{\bar{Q}}_{\omega}=\mathbf{I}_{3}$ & $\mathcal{\bar{Q}}_{b}=\mathbf{I}_{3}$\tabularnewline
		\hline 
		Case 2 & $\mathcal{\bar{Q}}_{v\left(i\right)}=0.1\mathbf{I}_{3}$ & $\mathcal{\bar{Q}}_{\omega}=10\mathbf{I}_{3}$ & $\mathcal{\bar{Q}}_{b}=10\mathbf{I}_{3}$\tabularnewline
		\hline 
		Case 3 & $\mathcal{\bar{Q}}_{v\left(i\right)}=0.01\mathbf{I}_{3}$ & $\mathcal{\bar{Q}}_{\omega}=100\mathbf{I}_{3}$ & $\mathcal{\bar{Q}}_{b}=100\mathbf{I}_{3}$\tabularnewline
		\hline 
	\end{tabular} 
\end{table}
In this Subsection, the color notation is as follows: black solid
and dashed lines represent the performance of the MEKF, blue center
line refers to the proposed direct attitude filter while orange and
purple depict the prescribed performance response. It can be noticed
in the upper portion of Figure \ref{fig:SO3PPF_9} that cases 1 and
2 show slower transient performance with less oscillatory behavior
in the steady-state. In contrast, the lower portion of Figure \ref{fig:SO3PPF_9}
illustrates that case 3 results in faster transient performance with
higher levels of oscillation in the steady-state. As such, a desired
measure of transient and stead-state error cannot be guaranteed in
case of MEKF. 
The direct attitude filter, on the other side, follows
the dynamically reducing boundaries achieving a desired level of prescribed
performance set by the user. 

The simulation results presented in this section validate the stable
performance and robustness of the two proposed filters against uncertain
measurements and large initialized errors. The two filters comply
with the constraints imposed by the user indicating guaranteed prescribed
performance measures in transient as well as steady-state performance.
This remarkable advantage was not offered in other nonlinear deterministic
attitude filters such as \cite{mahony2005complementary,mahony2008nonlinear,hamel2006attitude,hashim2018SO3Stochastic,grip2012attitude,hashim2018Conf1,lee2012exponential,zlotnik2017nonlinear}
as well as Gaussian attitude filters such as \cite{lefferts1982kalman,choukroun2006novel,markley2003attitude}.
Semi-direct attitude filter with prescribed performance requires attitude
reconstruction, for instance in our case we employed SVD \nameref{sec:SO3PPF_AppendixB},
to obtain $\tilde{R}=R_{y}^{\top}\hat{R}$. This adds complexity,
and therefore the semi-direct attitude filter requires more computational
power in comparison with direct attitude filter with prescribed performance.
However, both proposed filters showed remarkable convergence as detailed
in Table \ref{tab:SO3PPF_1}.

\section{Conclusion\label{sec:SO3PPF_Conclusion}}

In this paper, two nonlinear attitude filters with prescribed performance
characteristics have been considered. The filters are evolved directly
on $\mathbb{SO}\left(3\right)$. Attitude error has been defined in
terms of normalized Euclidean distance such that innovation term has
been selected to ensure predefined characteristics of transient and
steady-state performance. Consequently, the proposed filters achieve
superior convergence properties with transient error being less than
a predefined dynamic decreasing constrained function and steady-state
error being confined by a known lower bound. The constrained error
is transformed to its unconstrained form which is sufficient to solve
the attitude problem in prescribed performance sense. The filters
are deterministic while the stability analysis ensure boundedness
of all closed loop signals with asymptotic convergence of the normalized
Euclidean distance of attitude error to the origin. Simulation example
illustrated the robustness of the proposed filters in their response
to the predefined constraints in case when high level of uncertainties
is present in the measurements and a large initial attitude error
is observed.

\section*{Appendix A \label{sec:SO3PPF_AppendixA} }
\begin{center}
	\textbf{\large{}Proof of Lemma \ref{Lemm:SO3PPF_1}}{\large\par}
	\par\end{center}

Let the attitude be represented by $R\in\mathbb{SO}\left(3\right)$.
The attitude could be obtained knowing Rodriguez parameters vector
$\rho\in\mathbb{R}^{3}$. It is worth mentioning that Rodriguez parameters
vector $\rho$ is used for the sake of proving the results in Lemma
\ref{Lemm:SO3PPF_1}. The related map from vector form to $\mathbb{SO}\left(3\right)$
is governed by $\mathcal{R}_{\rho}:\mathbb{R}^{3}\rightarrow\mathbb{SO}\left(3\right)$
\cite{shuster1993survey} such that 
\begin{align}
\mathcal{R}_{\rho}\left(\rho\right)= & \frac{1}{1+||\rho||^{2}}\left(\left(1-||\rho||^{2}\right)\mathbf{I}_{3}+2\rho\rho^{\top}+2\left[\rho\right]_{\times}\right)\label{eq:SO3PPF_SO3_Rodr}
\end{align}
with direct substitution of \eqref{eq:SO3PPF_SO3_Rodr} in \eqref{eq:SO3PPF_Ecul_Dist}
one obtains 
\begin{equation}
||R||_{I}=\frac{||\rho||^{2}}{1+||\rho||^{2}}\label{eq:SO3PPF_TR2}
\end{equation}
Likewise, the anti-symmetric projection operator of attitude $R$
in \eqref{eq:SO3PPF_SO3_Rodr} for $\mathcal{R}_{\rho}=\mathcal{R}_{\rho}\left(\rho\right)$
can be defined as 
\begin{align*}
\boldsymbol{\mathcal{P}}_{a}\left(R\right)=\frac{1}{2}\left(\mathcal{R}_{\rho}-\mathcal{R}_{\rho}^{\top}\right)= & 2\frac{1}{1+||\rho||^{2}}\left[\rho\right]_{\times}
\end{align*}
and the vex operator 
\begin{equation}
\mathbf{vex}\left(\boldsymbol{\mathcal{P}}_{a}\left(R\right)\right)=2\frac{\rho}{1+||\rho||^{2}}\label{eq:SO3PPF_VEX_Pa}
\end{equation}
From \eqref{eq:SO3PPF_TR2} one can show that 
\begin{equation}
\left(1-||R||_{I}\right)||R||_{I}=\frac{||\rho||^{2}}{\left(1+||\rho||^{2}\right)^{2}}\label{eq:SO3PPF_append1}
\end{equation}
and from \eqref{eq:SO3PPF_VEX_Pa} one has 
\begin{equation}
||\mathbf{vex}\left(\boldsymbol{\mathcal{P}}_{a}\left(R\right)\right)||^{2}=4\frac{||\rho||^{2}}{\left(1+||\rho||^{2}\right)^{2}}\label{eq:SO3PPF_append2}
\end{equation}
Thus, \eqref{eq:SO3PPF_append1} and \eqref{eq:SO3PPF_append2} prove
\eqref{eq:SO3PPF_lemm1_1} in Lemma \ref{Lemm:SO3PPF_1}. From Subsection
\ref{subsec:SO3PPF_Explicit-Filter} $\sum_{i=1}^{n}s_{i}=3$ which
implies that ${\rm Tr}\left\{ M^{\mathcal{B}}\right\} =3$ and the
normalized Euclidean distance of $M^{\mathcal{B}}R$ is $||M^{\mathcal{B}}R||_{I}=\frac{1}{4}{\rm Tr}\left\{ M^{\mathcal{B}}\left(\mathbf{I}_{3}-R\right)\right\} $.
According to angle-axis parameterization in \eqref{eq:SO3PPF_att_ang},
one obtains

\begin{align}
||M^{\mathcal{B}}R||_{I} & =\frac{1}{4}{\rm Tr}\left\{ -M^{\mathcal{B}}\left(\sin(\theta)\left[u\right]_{\times}+\left(1-\cos(\theta)\right)\left[u\right]_{\times}^{2}\right)\right\} \nonumber \\
& =-\frac{1}{4}{\rm Tr}\left\{ M^{\mathcal{B}}\left(1-\cos(\theta)\right)\left[u\right]_{\times}^{2}\right\} \label{eq:SO3PPF_append3}
\end{align}
where ${\rm Tr}\left\{ M^{\mathcal{B}}\left[u\right]_{\times}\right\} =0$
as given in identity \eqref{eq:SO3PPF_Identity6}. One has \cite{murray1994mathematical}
\begin{equation}
||R||_{I}=\frac{1}{4}{\rm Tr}\left\{ \mathbf{I}_{3}-R\right\} =\frac{1}{2}\left(1-{\rm cos}\left(\theta\right)\right)={\rm sin}^{2}\left(\frac{\theta}{2}\right)\label{eq:SO3PPF_append4}
\end{equation}
and the Rodriguez parameters vector in terms of angle-axis parameterization
is \cite{shuster1993survey} 
\[
u={\rm cot}\left(\frac{\theta}{2}\right)\rho
\]
From identity \eqref{eq:SO3PPF_Identity3} $\left[u\right]_{\times}^{2}=-||u||^{2}\mathbf{I}_{3}+uu^{\top}$,
the expression in \eqref{eq:SO3PPF_append3} becomes 
\begin{align*}
||M^{\mathcal{B}}R||_{I} & =\frac{1}{2}||R||_{I}u^{\top}\bar{\mathbf{M}}^{\mathcal{B}}u\\
& =\frac{1}{2}||R||_{I}{\rm cot}^{2}\left(\frac{\theta}{2}\right)\rho^{\top}\bar{\mathbf{M}}^{\mathcal{B}}\rho
\end{align*}
From \eqref{eq:SO3PPF_append4}, one can find ${\rm cos}^{2}\left(\frac{\theta}{2}\right)=1-||R||_{I}$
which means 
\[
{\rm tan}^{2}\left(\frac{\theta}{2}\right)=\frac{||R||_{I}}{1-||R||_{I}}
\]
Consequently, the normalized Euclidean distance is defined in the
sense of Rodriguez parameters vector as 
\begin{align}
||M^{\mathcal{B}}R||_{I} & =\frac{1}{2}\left(1-||R||_{I}\right)\rho^{\top}\bar{\mathbf{M}}^{\mathcal{B}}\rho\nonumber \\
& =\frac{1}{2}\frac{\rho^{\top}\bar{\mathbf{M}}^{\mathcal{B}}\rho}{1+||\rho||^{2}}\label{eq:SO3PPF_append_MBR_I}
\end{align}
The anti-symmetric projection operator in terms of Rodriguez parameters
vector with aid of identity \eqref{eq:SO3PPF_Identity1} and \eqref{eq:SO3PPF_Identity4}
can be defined as 
\begin{align*}
\boldsymbol{\mathcal{P}}_{a}\left(M^{\mathcal{B}}R\right)= & \frac{M^{\mathcal{B}}\rho\rho^{\top}-\rho\rho^{\top}M^{\mathcal{B}}+M^{\mathcal{B}}\left[\rho\right]_{\times}+\left[\rho\right]_{\times}M^{\mathcal{B}}}{1+||\rho||^{2}}\\
= & \frac{\left[\left({\rm Tr}\left\{ M^{\mathcal{B}}\right\} \mathbf{I}_{3}-M^{\mathcal{B}}+\left[\rho\right]_{\times}M^{\mathcal{B}}\right)\rho\right]_{\times}}{1+||\rho||^{2}}
\end{align*}
It follows that the vex operator of the above expression is 
\begin{align}
\mathcal{\mathbf{vex}}\left(\boldsymbol{\mathcal{P}}_{a}\left(M^{\mathcal{B}}R\right)\right) & =\frac{1}{1+||\rho||^{2}}\left(\mathbf{I}_{3}-\left[\rho\right]_{\times}\right)\bar{\mathbf{M}}^{\mathcal{B}}\rho\label{eq:SO3PPF_append_MBR_VEX}
\end{align}
The 2-norm of \eqref{eq:SO3PPF_append_MBR_VEX} can be obtained by
\begin{align*}
||\mathbf{vex}\left(\boldsymbol{\mathcal{P}}_{a}\left(M^{\mathcal{B}}R\right)\right)||^{2} & =\frac{\rho^{\top}\bar{\mathbf{M}}^{\mathcal{B}}\left(\mathbf{I}_{3}-\left[\rho\right]_{\times}^{2}\right)\bar{\mathbf{M}}^{\mathcal{B}}\rho}{\left(1+||\rho||^{2}\right)^{2}}
\end{align*}
with the aid of identity \eqref{eq:SO3PPF_Identity3} $\left[\rho\right]_{\times}^{2}=-||\rho||^{2}\mathbf{I}_{3}+\rho\rho^{\top}$,
one obtains 
\begin{align}
||\mathbf{vex}\left(\boldsymbol{\mathcal{P}}_{a}\left(M^{\mathcal{B}}R\right)\right)||^{2} & =\frac{\rho^{\top}\bar{\mathbf{M}}^{\mathcal{B}}\left(\mathbf{I}_{3}-\left[\rho\right]_{\times}^{2}\right)\bar{\mathbf{M}}^{\mathcal{B}}\rho}{\left(1+||\rho||^{2}\right)^{2}}\nonumber \\
& =\frac{\rho^{\top}\left(\bar{\mathbf{M}}^{\mathcal{B}}\right)^{2}\rho}{1+||\rho||^{2}}-\frac{\left(\rho^{\top}\bar{\mathbf{M}}^{\mathcal{B}}\rho\right)^{2}}{\left(1+||\rho||^{2}\right)^{2}}\nonumber \\
& \geq\underline{\lambda}\left(1-\frac{\left\Vert \rho\right\Vert ^{2}}{1+\left\Vert \rho\right\Vert ^{2}}\right)\frac{\rho^{\top}\bar{\mathbf{M}}^{\mathcal{B}}\rho}{1+||\rho||^{2}}\label{eq:SO3PPF_append_VEX_ineq}
\end{align}
where $\underline{\lambda}=\underline{\lambda}\left(\bar{\mathbf{M}}^{\mathcal{B}}\right)$
is the minimum singular value of $\bar{\mathbf{M}}^{\mathcal{B}}$
and $\left\Vert R\right\Vert _{I}=\frac{\left\Vert \rho\right\Vert ^{2}}{1+\left\Vert \rho\right\Vert ^{2}}$
as defined in \eqref{eq:SO3PPF_TR2}. It can be found that 
\begin{align}
1-\left\Vert R\right\Vert _{I} & ={\rm Tr}\left\{ \frac{1}{12}\mathbf{I}_{3}+\frac{1}{4}R\right\} \nonumber \\
& ={\rm Tr}\left\{ \frac{1}{12}\mathbf{I}_{3}+\frac{1}{4}\left(M^{\mathcal{B}}\right)^{-1}M^{\mathcal{B}}R\right\} \label{eq:SO3PPF_append_rho2}
\end{align}
Therefore, from \eqref{eq:SO3PPF_append_VEX_ineq}, and \eqref{eq:SO3PPF_append_rho2}
the following inequality holds 
\begin{align*}
& ||\mathbf{vex}\left(\boldsymbol{\mathcal{P}}_{a}\left(M^{\mathcal{B}}R\right)\right)||^{2}\\
& \hspace{5em}\geq\frac{\underline{\lambda}}{2}\left(1+{\rm Tr}\left\{ \left(M^{\mathcal{B}}\right)^{-1}M^{\mathcal{B}}R\right\} \right)\left\Vert M^{\mathcal{B}}R\right\Vert _{I}
\end{align*}
which proves \eqref{eq:SO3PPF_lemm1_3} in Lemma \ref{Lemm:SO3PPF_1}.

\section*{Appendix B \label{sec:SO3PPF_AppendixB} }
\begin{center}
	\textbf{\large{}An Overview on SVD in }{\large{}\cite{markley1988attitude} }{\large\par}
	\par\end{center}

Let the true attitude be $R\in\mathbb{SO}\left(3\right)$. A set of
vectors presented in \eqref{eq:SO3PPF_Vector_norm} could be utilized
to reconstruct the attitude. Define $s_{i}$ as the confidence level
the $i$th measurement and for $n$ measurements. Define $s_{i}=\frac{s_{i}}{\sum_{i=1}^{n}s_{i}}$.
Accordingly, the corrupted reconstructed attitude $R_{y}$ is given
by
\[
\begin{cases}
\mathcal{J}\left(R\right) & =1-\sum_{i=1}^{n}s_{i}\left(\upsilon_{i}^{\mathcal{B}}\right)^{\top}R^{\top}\upsilon_{i}^{\mathcal{I}}\\
& =1-{\rm Tr}\left\{ R^{\top}B^{\top}\right\} \\
B & =\sum_{i=1}^{n}s_{i}\upsilon_{i}^{\mathcal{B}}\left(\upsilon_{i}^{\mathcal{I}}\right)^{\top}=USV^{\top}\\
U_{+} & =U\left[\begin{array}{ccc}
1 & 0 & 0\\
0 & 1 & 0\\
0 & 0 & {\rm det}\left(U\right)
\end{array}\right]\\
V_{+} & =V\left[\begin{array}{ccc}
1 & 0 & 0\\
0 & 1 & 0\\
0 & 0 & {\rm det}\left(V\right)
\end{array}\right]\\
R_{y} & =V_{+}U_{+}^{\top}
\end{cases}
\]
For more details consult \cite{markley1988attitude} or see the Appendix
in \cite{hashim2018SO3Stochastic}.

\section*{Appendix C \label{sec:SO3PPF_AppendixC} }
\begin{center}
	\textbf{\large{}An Overview of MEKF in }\cite{markley2003attitude,zamani2015nonlinear}{\large{} }{\large\par}
	\par\end{center}

The unit-quaternion vector $Q=\left[q_{0},q^{\top}\right]^{\top}\in\mathbb{S}^{3}$
is composed of a scalar component $q_{0}\in\mathbb{R}$ and a vector
$q\in\mathbb{R}^{3}$ defined by
\[
\mathbb{S}^{3}=\left\{ \left.Q\in\mathbb{R}^{4}\right|\left\Vert Q\right\Vert =1\right\} 
\]
The structure of MEKF is as follows
\begin{align}
\Psi\left(x\right) & =\left[\begin{array}{cc}
0 & -x^{\top}\\
x & -\left[x\right]_{\times}
\end{array}\right]\in\mathbb{R}^{4\times4},\hspace{1em}x\in\mathbb{R}^{3\times1}\nonumber \\
\dot{\hat{Q}} & =\frac{1}{2}\Psi\left(\Omega_{m}-\hat{b}+P_{a}W\right)\hat{Q}\nonumber \\
W & =\sum_{i=1}^{n}\hat{\upsilon}_{i}^{\mathcal{B}}\times\mathcal{\bar{Q}}_{v\left(i\right)}^{-1}\left(\hat{\upsilon}_{i}^{\mathcal{B}}-\upsilon_{i}^{\mathcal{B}}\right)\nonumber \\
\left[\begin{array}{c}
0\\
W
\end{array}\right] & =\sum_{i=1}^{n}\left[\begin{array}{c}
0\\
\hat{\upsilon}_{i}^{\mathcal{B}}
\end{array}\right]\times\left[\begin{array}{cc}
0 & 0_{3\times1}^{\top}\\
0_{3\times1} & \mathcal{\bar{Q}}_{v\left(i\right)}^{-1}
\end{array}\right]\left[\begin{array}{c}
0\\
\hat{\upsilon}_{i}^{\mathcal{B}}-\upsilon_{i}^{\mathcal{B}}
\end{array}\right]\nonumber \\
\left[\begin{array}{c}
0\\
\hat{\upsilon}_{i}^{\mathcal{B}}
\end{array}\right] & =\hat{Q}^{-1}\odot\left[\begin{array}{c}
0\\
\upsilon_{i}^{\mathcal{I}}
\end{array}\right]\odot\hat{Q}\label{eq:MEKF}
\end{align}
where $\hat{Q}\in\mathbb{S}^{3}$ denotes the estimate of the true
unit-quaternion, $\odot$ denotes unit-quaternion multiplication,
$\upsilon_{i}^{\mathcal{B}},\upsilon_{i}^{\mathcal{I}}\in\mathbb{R}^{3}$
are defined in \eqref{eq:SO3PPF_Vect_True} and \eqref{eq:SO3PPF_Vector_norm},
and
\begin{equation}
\begin{cases}
\boldsymbol{\mathcal{P}}_{s}\left(A\right) & =\frac{1}{2}\left(A+A^{\top}\right),\hspace{1em}A\in\mathbb{R}^{3\times3}\\
\dot{\hat{b}} & =P_{c}^{\top}W\\
S & =\sum_{i=1}^{n}\left[\hat{\upsilon}_{i}^{\mathcal{B}}\right]_{\times}\mathcal{\bar{Q}}_{v\left(i\right)}^{-1}\left[\hat{\upsilon}_{i}^{\mathcal{B}}\right]_{\times}\\
\dot{P}_{a} & =\mathcal{\bar{Q}}_{\omega}+2\boldsymbol{\mathcal{P}}_{s}\left(P_{a}\left[\Omega_{m}-\hat{b}\right]_{\times}-P_{c}\right)\\
& \hspace{1em}-P_{a}SP_{a}\\
\dot{P}_{b} & =\mathcal{\bar{Q}}_{b}-P_{c}SP_{c}\\
\dot{P}_{c} & =-\left[\Omega_{m}-\hat{b}\right]_{\times}P_{c}-P_{a}SP_{c}-P_{b}
\end{cases}\label{eq:MEKF-1}
\end{equation}
with $\mathcal{\bar{Q}}_{v\left(i\right)},\mathcal{\bar{Q}}_{\omega},\mathcal{\bar{Q}}_{b}\in\mathbb{R}^{3\times3}$
being covariance matrices, for $i=1,2,\ldots,n$. The rest of the
notation is identical to the notation used in the filter design given
before Theorem \ref{thm:SO3PPF_1} and \ref{thm:SO3PPF_2}.

\section*{Acknowledgment}

The authors would like to thank \textbf{Maria Shaposhnikova} for proofreading the article.

 \bibliographystyle{IEEEtran}
\bibliography{bib_PPF_SO3}

\clearpage

\section*{AUTHOR INFORMATION}
\vspace{10pt}
{\bf Hashim A. Hashim} is a Ph.D. candidate and a Teaching and Research Assistant in Robotics and Control, Department of Electrical and Computer Engineering at the University of Western Ontario, ON, Canada.\\
His current research interests include stochastic and deterministic filters on SO(3) and SE(3), control of multi-agent systems, control applications and optimization techniques.\\
\underline{Contact Information}: \href{mailto:hmoham33@uwo.ca}{hmoham33@uwo.ca}.
\vspace{50pt}

{\bf Lyndon J. Brown} received the B.Sc. degree from the U. of Waterloo, Canada in 1988 and the M.Sc. and PhD. degrees from the University of Illinois, Urbana-Champaign in 1991 and 1996, respectively. He is an associate professor in the department of electrical and computer engineering at Western University, Canada. He worked in industry for Honeywell Aerospace Canada and E.I. DuPont de Nemours.\\
His current research includes the identification and control of predictable signals, biological control systems, welding control systems, attitude and pose estimation.
\vspace{50pt}

{\bf Kenneth McIsaac} received the B.Sc. degree from the University of Waterloo, Canada, in 1996, and the M.Sc. and Ph.D. degrees from the University of Pennsylvania, in 1998 and 2001, respectively. He is currently an Associate Professor and the Chair of Electrical and Computer Engineering with Western University, ON, Canada. \\
His current research interests include computer vision and signal processing, mostly in the context of using machine intelligence in robotics and assistive systems. Also, his research interests include attitude and pose estimation.

\end{document}